\documentclass[a4pitaper,11pt]{amsart}
\usepackage{amsmath,amsthm,amssymb}
\usepackage[mathscr]{eucal}
 \usepackage{cite}
 \usepackage{bm} 
\usepackage{upgreek}
\usepackage[bookmarks,bookmarksnumbered]{hyperref}
\usepackage{enumerate}

\numberwithin{equation}{section}

\newcommand{\car}{\curvearrowright}

\theoremstyle{plain}
\newtheorem{main}{Theorem}
\newtheorem{mcor}[main]{Corollary}

\newtheorem{theorem}{Theorem}[section]

\newtheorem{lemma}[theorem]{Lemma}
\newtheorem{proposition}[theorem]{Proposition}
\newtheorem{corollary}[theorem]{Corollary}
\theoremstyle{definition}

\newtheorem{definition}[theorem]{Definition}
\newtheorem{example}[theorem]{Example}
\newtheorem{notation}[theorem]{Notation}
\newtheorem{remark}[theorem]{Remark}

\newtheorem{assumption}[theorem]{Assumption}
\begin{document}


\title[Measure equivalence rigidity  via s-malleable deformations]
{Measure equivalence rigidity  via s-malleable deformations}

\author{Daniel Drimbe}
\address{Department of Mathematics, KU Leuven, Celestijnenlaan 200b, B-3001 Leuven, Belgium}
\email{daniel.drimbe@kuleuven.be}
\thanks {The author holds the postdoctoral fellowship fundamental research 12T5221N of the Research Foundation Flanders.}

\begin{abstract} 
We single out a large class of groups $\bm{\mathscr{M}}$  for which the following unique prime factorization result holds: if $\Gamma_1,\dots,\Gamma_n\in \bm{\mathscr{M}}$ and $\Gamma_1\times\dots\times\Gamma_n$ is measure equivalent to a product $\Lambda_1\times\dots\times\Lambda_m$ of infinite icc groups, then $n \ge m$, and if $n = m$ then, after permutation of the indices, $\Gamma_i$ is measure equivalent to $\Lambda_i$, for all $1\leq i\leq n$. This provides an analogue of Monod and Shalom's theorem \cite{MS02} for groups that belong to $\bm{\mathscr{M}}$. Class $\bm{\mathscr{M}}$ is constructed using groups whose  von Neumann algebras admit an s-malleable deformation in the sense of Sorin Popa and it contains all icc non-amenable groups $\Gamma$ for which either (i) $\Gamma$ is an arbitrary wreath product group with amenable base or (ii) $\Gamma$ admits an unbounded 1-cocycle into its left regular representation. Consequently, we derive several orbit equivalence rigidity results for  actions of product groups that belong to $\bm{\mathscr{M}}$. 
Finally, for groups $\Gamma$ satisfying condition (ii), we show that all embeddings of group von Neumann algebras of non-amenable inner amenable  groups into $L(\Gamma)$ are ``rigid". In particular, we provide an alternative solution to a question of Popa that was recently answered in \cite{DKEP22}.

\end{abstract}

\maketitle

\section{Introduction}

Classifying countable groups up to {\it measure equivalence} is a central topic in measured group theory that has witnessed an explosion of activity for the last 25 years, see the surveys  \cite{Sh04,Fu09,Ga10} and the introduction of \cite{HH21}. The notion of measure equivalence has been introduced by Gromov \cite{Gr91} as a measurable analogue to the
geometric notion of quasi-isometry between finitely generated groups. Specifically, two countable groups $\Gamma$ and $\Lambda$ are called {\it measure equivalent} if there exist commuting free measure preserving actions of $\Gamma$ and $\Lambda$ on a standard measure space $(\Omega,m)$ such that the actions of  $\Gamma$ and $\Lambda$ on $(\Omega,m)$ each admit a finite measure fundamental domain. Natural examples of measure equivalent groups are two lattices in a locally compact second countable group. 

Measure equivalence can be studied through the lenses of orbit equivalence due to the fundamental result that two countable groups  are measure equivalent if and only if they admit free ergodic probability measure preserving (pmp) actions that are stably orbit equivalent \cite{Fu99}.
Recall that two pmp actions $\Gamma\car (X,\mu)$ and $\Lambda\car (Y,\nu)$ are called {\it stably orbit equivalent} if there exist non-null subsets $A\subset X, B\subset Y$ and a measure space isomorphism $\theta: A\to B$ such that $\theta(\Gamma x \cap A)= \Lambda \theta(x)\cap B$, for almost every $x\in A$. If $\mu(A)=\nu(B)=1$, then the two actions are called {\it  orbit equivalent} (OE). 

The celebrated work of Ornstein and Weiss \cite{OW80} (see also \cite{Dy59, CFW81}) proves that any two free ergodic pmp actions of infinite amenable groups are OE, and consequently, any two infinite amenable groups are measure equivalent. In sharp contrast, classifying non-amenable groups up to measure equivalence is a much more challenging task and it reveals a very strong rigidity phenomenon. By building on  Zimmer's work \cite{Zi84}, Furman showed that any countable group which is measure equivalent to a lattice in a higher rank simple Lie group is essentially a lattice in the same Lie group \cite{Fu99}. Then Kida showed that most mapping class groups ${\rm Mod}(S)$ are {\it measure equivalent superrigid} which means that any countable group that is measure equivalent  to ${\rm Mod}(S)$, must be virtually isomorphic to it \cite{Ki06}.
Subsequently, other such measure equivalent superrigid groups have been found and we refer the reader to  the introduction of \cite{HH22} for more details.



There have been discovered several other remarkable instances where various  aspects of the group $\Gamma$ can be recovered from its measure equivalence class or certain properties of the group action $\Gamma\car (X,\mu)$ are remembered by its associated orbit equivalence relation.
 We only highlight the following developments in this direction and refer the reader to the surveys  \cite{Sh04,Fu09} for more information. Gaboriau used the notion of cost to show that the rank of a free group $\mathbb F_n$ is an invariant of the  orbit equivalence relation of any of its free, ergodic, pmp actions \cite{Ga99}. Then  his discovery that measure equivalent groups have proportional $\ell^2$-Betti numbers \cite{Ga02} led to significant new progress in the classification problem of pmp actions up to OE, see the survey \cite{Ga10}.
Using a completely different conceptual framework, Popa's deformation rigidity/theory \cite{Po07} led to an unprecedented development in the theory of von Neumann algebras and provided many other spectacular rigidity results in orbit equivalence, see the surveys \cite{Va10a,Io12b,Io17}.

In their breakthrough work \cite{MS02}, Monod and Shalom employed techniques from bounded cohomology theory to obtain a series  of OE rigidity results, including the following unique prime factorization result: if $\Gamma_1\times\dots\times\Gamma_n$ is a product of non-elementary torsion-free hyperbolic groups (more generally, of groups belonging to class $\mathcal C_{\rm reg}$, see \cite[Notation 1.2]{MS02}) that  is measure equivalent to a product $\Lambda_1\times\dots\times\Lambda_m$ of torsion-free groups, then $n \ge m$, and if $n = m$ then, after permutation of the indices, $\Gamma_i$ is measure equivalent to $\Lambda_i$, for all $1\leq i\leq n$. By building upon  C$^*$-algebraic methods from \cite{Oz03,BO08}, the above unique prime factorization result has been extended by  Sako \cite{Sa09}  to products of non-amenable bi-exact groups (see also \cite{CS11}). 


In our first main result of the paper,
we use the powerful framework of Popa's deformation/rigidity theory to establish a general analogue of Monod and Shalom's unique prime factorization theorem, which applies in particular to product of groups with positive first $\ell^2$-Betti number. More generally, we obtain such a result for product of groups for which their von Neumann algebras belong to a certain class $\bm{\mathscr{M}}$ of II$_1$ factors that admit an s-malleable deformation in the sense of Popa \cite{Po01,Po03} (see Definition \ref{def:malleable}).
For simplicity, we say that a countable group $\Gamma$ belongs to $\bm{\mathscr{M}}$ if its associated von Neumann algebra $L(\Gamma)$ belongs to $\bm{\mathscr{M}}$.
We refer the reader to Definition \ref{definition.classM} for the description of class $\bm{\mathscr{M}}$ and to Example \ref{main.example} for more concrete examples of groups that belong to this class.

\begin{main}\label{B_0}
Let $\Gamma_1,\dots,\Gamma_n$ be groups that belong to $\bm{\mathscr{M}}$. If $\Gamma_1\times\dots\times\Gamma_n$ is measure equivalent to a product $\Lambda_1\times\dots\times\Lambda_m$ of infinite icc groups, then $n\ge m$, and if $n=m$, then after permutation of indices, $\Gamma_i$ is measure equivalent to $\Lambda_i$, for any $1\leq i\leq n$.
\end{main}

\begin{example}\label{main.example}
   A countable group $\Gamma$ belongs to $\bm{\mathscr{M}}$ whenever $\Gamma$ is a  non-amenable icc group that satisfies one of the following conditions (see Proposition \ref{remark.classM.classA}):
\begin{enumerate}
     \item $\Gamma=\Sigma\wr_{G/H} G$ is a generalized wreath product group with $\Sigma$ amenable, $G$ non-amenable and $H<G$ is an amenable almost malnormal subgroup.
    
    \item $\Gamma$ admits an unbounded cocycle for some mixing representation $\pi:\Gamma\to \mathcal O( H_{\mathbb R})$ such that $\pi$ is weakly contained in the left regular representation of $\Gamma$.
    
    \item $\Gamma=\Gamma_1*_\Sigma\Gamma_2$ is an amalgamated free product group satisfying $[\Gamma_1:\Sigma]\ge 2$ and $[\Gamma_2:\Sigma]\ge 3$, where $\Sigma<\Gamma$ is an amenable almost malnormal\footnote{A subgroup $H<G$ is called almost malnormal if $gHg^{-1}\cap H$ is finite for any $g\in G\setminus$H.} subgroup.

\end{enumerate}

\end{example}

We continue by making several remarks about Theorem \ref{B_0}. First, note  that the class $\mathcal C_{\rm reg}$ considered by Monod and Shalom in their unique prime factorization result \cite[Theorem 1.16]{MS02}  does not contain groups that have infinite amenable normal subgroups \cite[Corollary 1.19]{MS02}, and hence, the subclass of wreath product groups considered in Example \ref{main.example}(1) is disjoint from $\mathcal C_{\rm reg}$.  Moreover, Example \ref{main.example} provides a large class of groups that are not bi-exact \cite{Sa09} since any bi-exact group cannot contain an infinite subgroup with non-amenable centralizer.

Next, we contrast our result with the following corollary of Gaboriau's work \cite{Ga02}: if a product $\Gamma=\Gamma_1\times\dots\times\Gamma_n$ of $n$ groups with positive first $\ell^2$-Betti number is measure equivalent to a product $\Lambda=\Lambda_1\times\dots\times\Lambda_m$ of $m$ infinite groups, then $n\ge m$. Indeed, by \cite[Th\'{e}or\`{e}me 6.3.]{Ga02} we have that the $n^{\rm th}$ $\ell^2$-Betti number of $\Gamma$ vanishes if and only if the $n^{\rm th}$ $\ell^2$-Betti number of $\Lambda$ vanishes. On the other hand, 
the K\"{u}nneth formula \cite[Propri\'{e}t\'{e}s 1.5]{Ga02} implies that the $n^{\rm th}$ $\ell^2$-Betti number of $\Gamma$ is positive, while if $n<m$, then the $n^{\rm th}$ $\ell^2$-Betti number of $\Lambda$ equals to $0$.  Theorem \ref{A} strengthens this conclusion in two ways in the case $\Gamma,\Lambda$ are icc. First, if $n=m$ we are able to recover the measure equivalence class of each $\Gamma_i$. Secondly, since the groups with positive first $\ell^2$-Betti number are precisely the non-amenable groups that admit an unbounded cocycle into the left regular representation \cite{PT07}, Example \ref{main.example}(2) extends the previous result of Gaboriau to the larger class of groups that admit an unbounded cocycle for some mixing representation that is weakly contained in the left regular representation.

\begin{remark}
Popa's deformation/rigidity theory gave rise to a plethora of striking rigidity results for von Neumann algebras of wreath product groups. Popa's pioneering work \cite{Po03,Po04} allowed to distinguish between the group von Neumann algebras of $\mathbb Z/2\mathbb Z\wr \Gamma$, as $\Gamma$ is an infinite property (T) group, while Ioana, Popa and Vaes used  a wreath product construction to obtain the first class of groups that are entirely remembered by their von Neumann algebras \cite{IPV10}. Subsequently, several other rigidity results have been obtained for von Neumann algebras of wreath products  including primeness, relative solidity, and product rigidity, see \cite{Io06,Po06a,CI08,Io10,IPV10,SW11,CPS11,BV12,IM19, Dr20, CD-AD21}. Theorem \ref{B_0} establishes a new general rigidity result for wreath product groups  by showing that products of arbitrary non-amenable wreath product groups with amenable base satisfy an analogue of Monod and Shalom's unique prime factorization result.
\end{remark}








Theorem \ref{B_0} follows from the following more general result in which we classify  all tensor product decompositions of $L(\Lambda)$, whenever $\Lambda$ is an icc group that is measure equivalent to a finite product of groups that belong to  $\bm{\mathscr{M}}$.

\begin{main}\label{Theorem.ME}\label{A}
Let $\Gamma=\Gamma_1\times\dots\times\Gamma_n$ be a product of groups that belong to   $\bm{\mathscr{M}}$ and let $\Lambda$ be an icc group that is measure equivalent to  $\Gamma$. Assume $L(\Lambda)=P_1\bar\otimes \dots \bar\otimes P_m$ admits a tensor product decompositions into II$_1$ factors. Then $n\ge m$ and there exists a decomposition $\Lambda=\Lambda_1\times\dots\times\Lambda_m$ into infinite groups.

Moreover, there exist a partition $S_1\sqcup \dots \sqcup S_m=\{1,\dots,n\}$, a decomposition $L(\Lambda)=P_1^{t_1}\bar\otimes \dots\bar\otimes P_m^{t_m}$, for some $t_1,\dots,t_m>0$ with $t_1\dots t_m=1$, and a unitary $u\in L(\Lambda)$ such that for any $1\leq j\leq m$:
\begin{enumerate}
    \item $\times_{k\in S_j} \Gamma_k$ is measure equivalent to $\Lambda_j$;

    \item $P_j^{t_j}=uL(\Lambda_j)u^*$.
\end{enumerate}

In particular, if $n=m$, then after permutation of indices, $\Gamma_i$ is measure equivalent to $\Lambda_i$, for any $1\leq i\leq n$.
\end{main}

We note that Theorem \ref{A} provides a complement to \cite[Theorem C]{DHI16} where
such a classification result has been  obtained  by
Hoff, Ioana and the author in the case  the groups $\Gamma_i$ are hyperbolic.
Although the proof of Theorem \ref{A} is inspired by the strategy of the proof of 
\cite[Theorem C]{DHI16}, we implement quite differently some of the steps. In order to effectively work with groups  from $\bm{\mathscr{M}}$, which  are  defined by a property of their von Neumann algebras, we are making use in an essential way of newer techniques from \cite{BMO19,IM19,Dr20}. In particular,  our proof uses a relative version of the flip automorphism method introduced by Isono and Marrakchi in \cite{IM19}




Another application of Thereom \ref{A} is to the study of tensor product decompositions of von Neumann algebras by providing new classes of prime II$_1$ factors. Recall that a II$_1$ factor  is called {\it prime} if it does not admit a tensor product decomposition into II$_1$ factors. Popa discovered in \cite{Po83} the first examples of prime II$_1$ factors by showing that the free group factors $L(\mathbb F_S)$, with $S$ uncountable, are prime. Then Ge  showed in \cite{Ge96} that the free group factors $L(\mathbb F_n), 2\leq n\leq \infty$, are prime, thus providing the first examples of separable prime II$_1$ factors. Subsequently, a large number of prime II$_1$ factors have been discovered, see for instance the introduction of \cite{CDI22}.
As a corollary of Theorem \ref{A}, we obtain that if $\Gamma$ is a countable group that belongs to $\bm{\mathscr{M}}$ and $G=(\times_{i=1}^n \Gamma)\rtimes \mathbb Z/ n\mathbb Z$ is the semidirect product group of the natural translation action $\mathbb Z/n\mathbb Z \car \times_{i=1}^n \Gamma$, then $L(G)$ is a prime II$_1$ factor. In fact, a more general result holds and for properly formulating it, we make the following notation. Let $n$ be a positive integer, denote by $S_n$ the group of permutations of $\{1,\dots, n\}$ and consider the permutation action of $S_n$ on $\{1,\dots,n\}$. 
For any subset $J\subset \{1,\dots,n\}$ and subgroup $K<S_n$, we denote ${\rm Fix}_K(J)=\{g\in K\; | \; g\cdot j=j, \text{ for any } j\in J  \}$.

\begin{mcor}\label{Main.new.prime}
Let $\Gamma$ be a countable group that belongs to $\bm{\mathscr{M}}$. Let $n$ be a positive integer and let $K$ be any subgroup of $S_n$. Consider the permutation action $K\car \times_{i=1}^n\Gamma$ and denote $G=(\times_{i=1}^n\Gamma)\rtimes K$.

Then $L(G)$ is  a prime II$_1$ factor if and only if there exists no partition $J_1\sqcup J_2=\{1,\dots,n\}$ for which $K={\rm Fix}_K(J_1)\times {\rm Fix}_K(J_2)$.
\end{mcor}

Note that Corollary \ref{Main.new.prime} provides a large class of prime II$_1$ factors which admit finite index subfactors that are not prime. Additional such prime II$_1$ factors have been previously obtained in  \cite{DHI16,CD19} by replacing $\Gamma$ 
in the statement of Corollary \ref{Main.new.prime} by any non-elementary hyperbolic group, see also \cite[Section 5]{CDI22}.




    
    
    

    
    

We continue by discussing some OE-rigidity results for actions of product groups that belong to class $\bm{\mathscr{M}}$. Furman discovered in \cite{Fu99} the first class of group action $\Gamma\car (X,\mu)$ that are {\it OE-superrigid}, that is, any free, ergodic, pmp action that is OE to $\Gamma\car(X,\mu)$ must be virtually conjugate\footnote{Two pmp actions $\Gamma\car (X,\mu)$ and $\Lambda\car (Y,\nu)$ are {\it virtually conjugate} if there exist some finite normal subgroups $A<\Gamma$ and $B<\Lambda$ such that the associated actions $\Gamma/A\car X/A$ and $\Lambda/B\car Y/B$ are induced from conjugate actions.} to it. Subsequently, a large number of OE-superrigidity results have been obtained, see the introduction of \cite{DIP19}. 
By using part of the proof of Theorem \ref{A}  together with results from measured group theory \cite{HH21}, we derive the following  OE-superrigidity result within the class of mildly mixing actions.
Before stating the result,
 we recall  some notions. A pmp action $\Gamma_1\times \dots\times \Gamma_n\car (X,\mu)$ is called {\it irreducible} if its restriction to any subgroup $\Gamma_i$ is ergodic. A pmp action $\Lambda\car (Y,\nu)$ is called {\it mildly mixing} if whenever $A\subset Y$ is measurable subset satisfying $\liminf_{g\to\infty} \nu(gA\Delta A)=0$, then $\nu(A)\in\{0,1\}.$

\begin{main}\label{C}
Let $\Gamma=\Gamma_1\times\dots\times\Gamma_n$ be a product of $n\ge 2$ groups that belong to $\bm{\mathscr{M}}$.  Let $\Gamma\car(X,\mu)$ be a free, irreducible, pmp action that is orbit equivalent to a free, mildly mixing, pmp action $\Lambda\car (Y,\nu)$.

Then $\Gamma\car(X,\mu)$ and $\Lambda\car (Y,\nu)$ are virtually conjugate.

\end{main}

Note that this type of superrigidity has been obtained by Monod and Shalom  \cite[Theorem 1.9]{MS02} for groups $\Gamma_i$ that are torsion-free hyperbolic groups (more generally, groups that belong to $\mathcal C_{\rm reg}$). In Theorem \ref{C} we extend this result to groups from $\bm{\mathscr{M}}$ which are purely defined by a property of their von Neumann algebra.

Finally, in our last part of the paper we discuss some structural results for II$_1$ factors that belong to a subclass of $\bm{\mathscr{M}}$.
We say that a non-amenable tracial von Neumann algebra $M$ belongs to class $\bm{\mathscr{M}_0}$ if there exists an s-malleable deformation $(\tilde M, (\alpha_t)_{t\in\mathbb R})$ of $M$  satisfying:
\begin{itemize}
        \item   $L^2(\tilde M)\ominus L^2(M)$ is a mixing $M$-$M$-bimodule relative to $\mathbb C1$.
    
        \item $L^2(\tilde M)\ominus L^2(M)$ is weakly contained in the coarse bimodule $L^2(M)\otimes L^2(M)$ as $M$-$M$-bimodules.    
\end{itemize}

We refer the reader to Sections \ref{subsection.bimodules} and \ref{subsection.malleable} for the terminology used in defining the class $\bm{\mathscr{M}_0}$ and we note that any II$_1$ factor from $\bm{\mathscr{M}_0}$ belongs to $\bm{\mathscr{M}}$, see Definition \ref{definition.classM}. 

In Theorem \ref{D} below we show that all embeddings of group von Neumann algebras of non-amenable inner amenable  groups in any II$_1$ factor that belongs to $\bm{\mathscr{M}_0}$ are {\it rigid}.
A countable group $\Gamma$ is {\it inner amenable } if there exists an atomless mean on $\Gamma$ which is invariant by the action of $\Gamma$ on itself by conjugation. Effros made in \cite{Ef75} a connection of this group theoretic notion to von Neumann algebras by showing that an icc group $\Gamma$ is inner amenable  whenever its group von Neumann algebra has property Gamma. The converse is false as was shown by Vaes \cite{Va09}.

\begin{main}\label{main.theorem.inneramenable}\label{D}
Let $M$ be a von Neumann algebra in $\bm{\mathscr{M}_0}$ and let $(\tilde M, (\alpha_t)_{t\in\mathbb R})$ be the associated s-malleable deformation of $M$. Let $\Gamma$ be a non-amenable inner amenable  group satisfying $L(\Gamma)\subset M$. 

Then $L(\Gamma)$ is $\alpha$-rigid, i.e. $\alpha_t\to {\rm id}$ uniformly on the unit ball of $L(\Gamma)$.
\end{main}

Note that von Neumann algebras with property Gamma exhibit strong structural results (see, for instance, \cite{Pe06, HU15, IS18}) that are enough for obtaining various rigidity results via 
Popa’s deformation/rigidity theory. In order to work with the more general class of inner amenable groups, we use an idea from  \cite[Theorem 11]{TD14} on how to use Popa's spectral gap principle. An additional obstacle that arises here is the fact that $E_{L(\Gamma)}(\alpha_t(u_g))$ is not necessarily a scalar multiple of $u_g$, where $g\in\Gamma$; here, we denoted by $\{u_g\}_{g\in\Gamma}$ the canonical unitaries that generate $L(\Gamma)$ and by $E_{L(\Gamma)}:\tilde M\to L(\Gamma)$ the canonical conditional expectation.
We overcome this difficulty by using an augmentation technique based on the commultiplication map associated to $L(\Gamma)$ \cite{PV09}.

We continue by discussing several applications of Theorem \ref{D}.
Chifan and Sinclair proved in \cite{CS11} that any countable group $\Gamma$ for which $\beta^{(2)}_1(\Gamma)>0$ is not inner amenable. Theorem \ref{D}  recovers and strengthens this fact in the following way. While it is unknown that the non-vanishing of the first $\ell^2$-Betti number is a group von Neumann algebra invariant, we derive from Theorem \ref{D} that any group that has isomorphic von Neumann algebra to $L(\Gamma)$ is not inner amenable  as well.

\begin{mcor}\label{mcor.inner.amenable}
Let $\Gamma$ be any countable group for which $\beta^{(2)}_1(\Gamma)>0$. If $\Lambda$ is any countable group for which $L(\Gamma)\cong L(\Lambda)$, then $\Lambda$ is not inner amenable.

\end{mcor}

To put Theorem \ref{D} into a better perspective, we note that
it provides an alternative solution to a question of Popa. Since any non-amenable property Gamma von Neumann algebra cannot embed into the free group factor $L(\mathbb F_n)$ \cite{Oz03},  Popa asked in \cite{Po21} if it still true that the group von Neumann algebra of a nonamenable inner amenable  group cannot embed into $L(\mathbb F_n)$. Recently, inspired by the notion of properly proximal groups \cite{BIP18} (see also \cite{IPR19}), Ding, Kunnawalkam Elayavalli and Peterson developed subtle boundary techniques to define a notion of  proper proximality  for tracial von Neumann algebras, and as a consequence, they answered Popa's question in a positive way. As a particular case of Theorem \ref{D}, we give a new proof for Popa's question by using methods from Popa's deformation/rigidity theory. 



Moreover, as a corollary  of Theorem \ref{D} we completely classify all embeddings of group von Neumann algebras $L(G)$ of non-amenable inner amenable  groups in any free product $M=M_1*M_2$ of tracial von Neumann algebras by showing that $L(G)\prec_{M} M_i$, for some $i$. Here, $\prec_M$ refers to Popa's intertwining-by-bimodules technique, see Section \ref{section.popa.intertwining}. Consequently, we obtain a new class of examples for which the Kurosh-type rigidity results discovered in \cite{Oz04} for free products von Neumann algebras hold. Namely, Ozawa
proved using C$^*$-algebraic techniques that if there is an isomorphism $\theta: M_1*\dots* M_m\to N_1*\dots* N_n$, where all von Neumann algebras $M_i$ and $N_j$ are non-amenable, semiexact, non-prime II$_1$ factors, then $m=n$, and after a permutation of indices, $\theta(M_i)$ is unitarily conjugate to $N_i$, for any $i\in\overline{1,n}$. By using Popa's deformation/rigidity theory, Ioana, Popa and Peterson obtained the previous Kurosh-type rigidity result for property (T) II$_1$ factor.
Shortly after, by developing a new approach rooted on closable derivations, Peterson unified and generalized these Kurosh-type rigidity  results  by covering $L^2$-rigid II$_1$ factors, which include all non-amenable non-prime, property (T) and property Gamma II$_1$ factors  \cite{Pe06}. By classifying certain amenable subalgebras of amalgamated free product von Neumann algebras, Ioana then extended the previous Kurosh-type rigidity result by covering non-amenable II$_1$ factors that admit a Cartan subalgebra  \cite{Io12a}.  
We also refer the reader to \cite{HU15} for certain Kurosh-type rigidity results for type III factors. 
As a corollary of Theorem \ref{D}, we  extend the previous  Kurosh-type rigidity results  to the class of II$_
1$ factors of non-amenable inner amenable  groups, see Corollary \ref{corollary.kurosh}.






{\bf Acknowledgements.} I am grateful to Adrian Ioana for sharing with me his original proof for \cite[Theorem 11]{TD14} which led to Theorem \ref{D}. I would like to thank Adrian Ioana and Stefaan Vaes for numerous comments and suggestions that helped improve the exposition of the paper. I would also like to thank   Srivatsav Kunnawalkam Elayavalli and  Changying Ding for  helpful comments.

\section{Preliminaries}

\subsection{Terminology}

Throughout the paper we consider {\it tracial von Neumann algebras} $(M,\tau)$, i.e. von Neumann algebras $M$ equipped with a faithful normal tracial state $\tau: M\to\mathbb C.$ This induces a norm on $M$ by the formula $\|x\|_2=\tau(x^*x)^{1/2},$ for any $x\in M$. We will always assume that $M$ is {\it separable}, i.e. the $\|\cdot\|_2$-completion of $M$ denoted by $L^2(M)$ is separable as a Hilbert space.
We denote by $\mathcal Z(M)$ the {\it center} of $M$ and by $\mathcal U(M)$ its {\it unitary group}. For two von Neumann subalgebras $P_1,P_2\subset M$, we denote by $P_1\vee P_2=W^*(P_1\cup P_2)$ the von Neumann algebra generated by $P_1$ and $P_2$. 

All inclusions $P\subset M$ of von Neumann algebras are assumed unital. We denote by $E_{P}:M\to P$ the unique $\tau$-preserving {\it conditional expectation} from $M$ onto $P$, by $e_P:L^2(M)\to L^2(P)$ the orthogonal projection onto $L^2(P)$ and by $\langle M,e_P\rangle$ the Jones' basic construction of $P\subset M$. We also denote by $P'\cap M=\{x\in M|xy=yx, \text{ for all } y\in P\}$ the {\it relative commutant} of $P$ in $M$ and by $\mathcal N_{M}(P)=\{u\in\mathcal U(M)|uPu^*=P\}$ the {\it normalizer} of $P$ in $M$.

The {\it amplification} of a II$_1$ factor $(M,\tau)$ by a number $t>0$ is defined to be $M^t=p(\mathbb B(\ell^2(\mathbb Z))\bar\otimes M)p$, for a projection $p\in \mathbb B(\ell^2(\mathbb Z))\bar\otimes M$ satisfying $($Tr$\otimes\tau)(p)=t$. Here Tr denotes the usual trace on $\mathbb B(\ell^2(\mathbb Z))$. Since $M$ is a II$_1$ factor, $M^t$ is well defined. Notice that if $M=P_1\bar\otimes P_2$, for some II$_1$ factors $P_1$ and $P_2$, then there is a natural isomorphism $M=P_1^t\bar\otimes P_2^{1/t}$, for any $t>0.$


Finally, for a positive integer $n$, we denote by $\overline{1,n}$ the set $\{1,\dots, n\}$. If $S\subset \overline{1,n}$ we denote its complement by $\widehat S=\overline{1,n}\setminus S$. In the case that $S=\{i\},$ we will simply write $\hat i$ instead of $\widehat {\{i\}}$. Also, given
any product group $\Gamma=\Gamma_1\times \dots\times \Gamma_n$, we will denote their subproduct supported on $S$ by $\Gamma_S=\times_{i\in S}\Gamma_i$.

\subsection {Intertwining-by-bimodules}\label{section.popa.intertwining} We next recall from  \cite [Theorem 2.1 and Corollary 2.3]{Po03} the powerful {\it intertwining-by-bimodules} technique of S. Popa.

\begin {theorem}[\!\!\cite{Po03}]\label{corner} Let $(M,\tau)$ be a tracial von Neumann algebra and $P\subset pMp, Q\subset qMq$ be von Neumann subalgebras. Let $\mathcal U\subset\mathcal U(P)$ be a subgroup such that $\mathcal U''=P$.

Then the following are equivalent:

\begin{enumerate}

\item There exist projections $p_0\in P, q_0\in Q$, a $*$-homomorphism $\theta:p_0Pp_0\rightarrow q_0Qq_0$  and a non-zero partial isometry $v\in q_0Mp_0$ such that $\theta(x)v=vx$, for all $x\in p_0Pp_0$.

\item There is no sequence $(u_n)_n\subset\mathcal U$ satisfying $\|E_Q(xu_ny)\|_2\rightarrow 0$, for all $x,y\in M$.

\end{enumerate}
\end{theorem}

If one of the equivalent conditions of Theorem \ref{corner} holds true, we write $P\prec_{M}Q$, and say that {\it a corner of $P$ embeds into $Q$ inside $M$.}
If $Pp'\prec_{M}Q$ for any non-zero projection $p'\in P'\cap pMp$, then we write $P\prec^{s}_{M}Q$.

\begin{lemma}\label{lemma.joint.control}
Let $\Lambda\car B$ be a trace preserving action and denote $M=B\rtimes\Lambda$. Let $p\in B$ be a non-zero projection and let $A\subset pBp$ be a von Neumann subalgebra such that $A'\cap pMp\subset A$. 

Let $\Lambda_0<\Lambda$ be a subgroup and $\mathcal G\subset\mathcal N_{pMp}(A)$ a group of unitaries.
If there is a projection $e\in \mathcal G'\cap pMp$ satisfying $\mathcal G''e\prec_M^s B\times\Lambda_0$, then there is a projection $f\in (A\cup\mathcal G)'\cap pMp$ with $ e\leq f$ satisfying $(A\cup\mathcal G)''f\prec_M^s B\rtimes\Lambda_0$. 
\end{lemma}


{\it Proof.} Throughout the proof we use the notation that a set $F\subset \Lambda$ is said to be small relative to $\{\Lambda_0\}$ if it is contained into a finite union of $s\Lambda_0 t$, where $s,t\in\Lambda$.
For any $F\subset\Lambda$, let $\mathcal H_F\subset L^2(M)$ be the $\|\cdot \|_2$-closed linear span of $\{B v_{\lambda}|\lambda\in F\}$ and denote by $P_F:L^2(M)\to \mathcal H_F$ the orthogonal projection onto $\mathcal H_F$. Let $\epsilon>0$ and denote $T=E_{A'\cap pMp}(e)$. Note that $T\in (A\cup\mathcal G)'\cap pMp\subset A$ and that $T$ belongs to the $\|\cdot\|_2$-closed convex hull of $\{aea^*|a\in \mathcal U(A)\}$. Thus, we can take $a_1,\dots,a_n\in\mathcal U(A)$ and $\alpha_1,\dots,\alpha_n\in [0,1]$ such that if we denote $T_0=\sum_{i=1}^n \alpha_i a_i ea_i^*$, then $\|T-T_0\|_2\leq \epsilon$.

Since $\mathcal G''e\prec_M^s B\times\Lambda_0$, it follows from \cite[Lemma 2.5]{Va10b} that there exists  $F\subset\Lambda$ that is small relative to $\{\Lambda_0\}$ such that $\|we-P_{F}(we)\|_2\leq \epsilon/n$, for all $w\in \mathcal G$. Hence, for all $a\in \mathcal U(A), w\in\mathcal G$, we have
\begin{equation}\label{x1}
    \|awT - \alpha_i\sum_{i=1}^n a (wa_iw^*)P_F(we)a_i^*  \|_2\leq \epsilon + \sum_{i=1}^n \alpha_i\| aw(a_iea_i^*) - a(wa_iw^*)P_F(we)  a_i^* \|_2\leq 2\epsilon.
\end{equation}

Since $A\subset pBp$, we have $a (wa_iw^*)P_F(we)a_i^* \in \mathcal H_F$, and hence, $\|awT-P_F(awT)\|_2\leq 2\epsilon$, for all $a\in\mathcal U(A),w\in\mathcal G$.
Therefore, there exists a sequence $\{F_n\}_{n\ge 1}$ of subsets of $\Lambda$ that are small relative to $\{\Lambda_0\}$ such that $\|awT-P_{F_n}(awT)\|_2\to 0$ uniformly in $a\in \mathcal U(A),w\in \mathcal G$. 

For every $\delta>0$ define the spectral projection $q_\delta=\chi_{(\delta,\infty)}(T)\in A$ and let $T_\delta\in A$ satisfying $TT_\delta=q_\delta.$ If we denote by $q_0$ the support projection of $T$, then $\|q_\delta-q_0\|_2\to 0$ as $\delta\to 0$. These altogether imply that $\|awq_\delta-P_{F_n}(awq_\delta)\|_2\to 0$ uniformly in $a\in \mathcal U(A),w\in \mathcal G$, and hence, $\|awq_0-P_{F_n}(awq_0)\|_2\to 0$ uniformly in $a\in \mathcal U(A),w\in \mathcal G$.
Finally, note that $q_0\in (A\cup\mathcal G)'\cap pMp$ and $q_0\ge e$. This concludes the proof. 
\hfill$\blacksquare$

The following proposition follows from \cite{BMO19,IM19} and it is essentially contained in the proof of \cite[Theorem 4.2]{Dr20}. We record it here for the convenience of the reader.

\begin{proposition}\label{proposition.full.sequence}
Let $M=P\bar\otimes Q$ be a tensor product of II$_1$ factors. Let $Q_n,n\ge 1$, be a decreasing sequence of von Neumann subalgebras such that $P\prec_M \bigvee_{n\ge 1}(Q_n'\cap M)$.

If $P$ does not have property Gamma, then there exists $m\ge 1$ such that $P\prec_{M} Q_m'\cap M.$
\end{proposition}

Recall that a II$_1$ factor $(M,\tau)$ has { property Gamma} if it admits a central sequence $(x_n)_n\subset{\mathcal U}(M)$ for which ${\rm inf}_{n} \|x_n- \tau(x_n)1\|_2>0$. 


\subsection{Bimodules}\label{subsection.bimodules}
Let $M, N$ be tracial von Neumann algebras. An $M$-$N$ bimodule $_M \mathcal H_N$ is a Hilbert space $\mathcal H$ together with with a $*$-homomorphism $\pi_{\mathcal H}: M\odot N^{\text{op}}\to\mathbb B(\mathcal H)$ that is normal on $M$ and $N^{\text{op}}$, where $M\odot N^{\text{op}}$ is the algebraic tensor product between $M$ and the opposite von Neumann algebra $N^{\text{op}}$ of $N$. Examples of bimodules include the trivial $M$-bimodule $_M L^2(M)_M$ and the coarse $M$-$N$-bimodule $_M L^2(M)\otimes L^2(N)_N.$
For two $M$-$N$-bimodules $_M \mathcal H_N$ and $_M \mathcal K_N$, we say that $_M \mathcal H_N$ is {\it weakly contained} in $_M \mathcal K_N$ if $\|\pi_{\mathcal H}(x)\|\leq \|\pi_{\mathcal K}(x)\|$, for any $x\in M\odot N^{\text{op}}$. 

Let $A\subset M$ be an inclusion of tracial von Neumann algebras and let $_M\mathcal H_M$ be an $M$-bimodule. We say that $_M\mathcal H_M$ is {\it mixing relative} to $A$ if for any sequence $(x_n)_n\subset (M)_1$ satisfying $\|E_{A}(xu_n y)\|_2\to 0$, for all $x,y\in M$, we have 
$$\lim_{n\to\infty} \sup_{y\in (M)_1}  \langle x_n\xi y,\eta \rangle, \text{  for all }\xi,\eta\in\mathcal H.$$

\subsection{Relative amenability}

A tracial von Neumann algebra $(M,\tau)$ is {\it amenable} if there is a positive linear functional $\Phi:\mathbb B(L^2(M))\to\mathbb C$ such that $\Phi_{|M}=\tau$ and $\Phi$ is $M$-{\it central}, meaning $\Phi(xT)=\Phi(Tx),$ for all $x\in M$ and $T\in \mathbb B(L^2(M))$. By Connes' classification of amenable factors \cite{Co76}, it follows that $M$ is amenable if and only if $M$ is approximately finite dimensional.
 
We continue by recalling the notion of relative amenability which is due to 
Ozawa and Popa \cite{OP07}. Fix a tracial von Neumann algebra $(M,\tau)$. Let $p\in M$ be a projection and $P\subset pMp,Q\subset M$ be von Neumann subalgebras. Following \cite[Definition 2.2]{OP07}, we say that $P$ is {\it amenable relative to $Q$ inside $M$} if there is a positive linear functional $\Phi:p\langle M,e_Q\rangle p\to\mathbb C$ such that $\Phi_{|pMp}=\tau$ and $\Phi$ is $P$-central. 
We say that $P$ is {\it strongly non-amenable relative to} $Q$ if $Pp'$ is non-amenable relative to $Q$ for any non-zero projection $p'\in P'\cap pMp$ (equivalently, for any non-zero projection $p'\in \mathcal{N}_M(P)'\cap pMp$ by \cite[Lemma 2.6]{DHI16}).

Note that if $P\subset pMp$ and $Q\subset M$ are tracial von Neumann algebras, then $P\subset pMp$ is amenable relative to $Q$ if and only if $_P L^2(pM)_M$ is weakly contained in $_P L^2(p\langle M,e_{Q} \rangle)_M$. We also recall that $_M L^2(\langle M,e_{Q} \rangle)_M\cong_M( L^2(M)\otimes_Q L^2(M))_M$. 
It is clear that $P$ is amenable relative to $\mathbb C 1$ inside $M$ if and only if $P$ is amenable.
The following lemma generalizes this fact and it is inspired by the proof of \cite[Lemma 5.6]{DHI16}. For completeness, we provide all the details.

\begin{lemma}\label{lemma.tensor}

Let $M_0$ and $M\subset \tilde M$ be some tracial von Neumann algebras, and $Q\subset q(M_0\bar\otimes M)q$ a von Neumann subalgebra. The following hold:

\begin{enumerate}

\item Assume that $Q\tilde z$ is amenable relative to $M_0$ inside $M_0\bar\otimes \tilde M$, for a non-zero projection $ \tilde z\in Q'\cap q(M_0\bar\otimes \tilde M)q$. Then $Qz$ is amenable relative to $M_0$ inside $M_0\bar\otimes M$, where   $z\in Q'\cap q(M_0\bar\otimes M)q $ is the support projection of $E_{\mathcal M}(\tilde z).$ 

\item If $Q\prec _{M_0\bar\otimes \tilde M} M_0,$ then $Q\prec_{M_0\bar\otimes M} M_0.$

\end{enumerate}

\end{lemma}

{\it Proof.} (1) Let $\mathcal M=M_0\bar\otimes M$ and $\tilde{\mathcal M}=M_0\bar\otimes\tilde M.$ The assumption implies that the bimodule $_{Q\tilde z} L^2(\tilde z\tilde{\mathcal M})_{\tilde{\mathcal M}}$ is weakly contained in $_{Q\tilde z} L^2(\tilde z \langle \tilde{\mathcal M},e_{M_0}\rangle)_{\tilde{\mathcal M}}$. If we denote by $z\in Q'\cap q(M_0\bar\otimes M)q)$ the support projection of $E_{\mathcal M}(\tilde z),$ we obtain that
\begin{equation}\label{eqq1}
_{Qz} L^2(z{\mathcal M})_{\mathcal M} \text{ is weakly contained in }_{Qz} L^2(z \langle \tilde{\mathcal M},e_{M_0}\rangle)_{{\mathcal M}}.    
\end{equation}

Notice that $_{\mathcal M} L^2(\langle \tilde{\mathcal M},e_{M_0}\rangle)_\mathcal M\cong$ $ _{\mathcal M}  L^2(\tilde {\mathcal M})\otimes L^2(\tilde M)  _{\mathcal M}$. Note also that $_{\mathcal M} L^2(\tilde{\mathcal M})_{M_0}$ is weakly contained in $_{\mathcal M} L^2({\mathcal M})_{M_0}$ and
 $_{\mathbb C} L^2(\tilde M)_M$ is weakly contained in $_{\mathbb C}L^2(M)_M.$ 

These altogether imply that  $_{Qz} L^2(\langle \tilde{\mathcal M},e_{M_0}\rangle)_\mathcal M$ is weakly contained in $_{Qz} L^2(\langle {\mathcal M},e_{M_0}\rangle)_\mathcal M$. Using \eqref{eqq1} we deduce that $_{Qz} L^2(z{\mathcal M})_{{\mathcal M}}$ is weakly contained in $_{Qz} L^2(z \langle {\mathcal M},e_{M_0}\rangle)_{{\mathcal M}}$, which shows that $Qz$ is amenable relative to $M_0$ inside $\mathcal M.$

(2) By assuming the contrary,  there exists a sequence $u_n\in \mathcal U (Q)$ such that 
\begin{equation}\label{larger}
\|E_{M_0}(xu_ny)\|_2\to 0, \text{ for all } x,y\in M_0\bar\otimes M.
\end{equation}
We want to show that \eqref{larger} holds for all $x,y\in M_0\bar\otimes \tilde M,$ which will contradict the assumption. Note that it is enough to consider $x=1$ and $y\in\tilde M$. In this case, by using \eqref{larger} we obtain 
$E_{M_0}(u_ny)=E_{M_0}(E_{M_0\bar\otimes M}(u_ny))=E_{M_0}(u_nE_{M_0\bar\otimes M}(y))$, which goes to $0$ in the $\|\cdot\|_2$-norm. This finishes the proof.
\hfill$\blacksquare$

\section{Malleable deformations for von Neumann algebras: class $\mathscr{M}$}

\subsection{Malleable deformations}\label{subsection.malleable}
 Popa introduced in  \cite{Po01, Po03} the notion of an s-malleable deformation of a von Neumann algebra. This notion has been successfully used in the framework of his deformation/rigidity theory and led to a plethora of remarkable results in the theory of von Neumann algebras, see
 the surveys \cite{Po07,Va10a, Io12b,Io17}. We also refer the reader to \cite{dSHHS20} for recent developments on s-malleable deformations.

\begin{definition}\label{def:malleable}
Let $(M,\tau)$ be a tracial von Neumannn algebra. A pair $(\tilde M, (\alpha_t)_{t\in\mathbb R})$ is called an {\it s-malleable deformation} of $M$ if the following conditions hold:
\begin{itemize}
    \item $(\tilde M,\tilde\tau)$ is a tracial von Neumann algebra such that $M\subset\tilde M$ and $\tau=\tilde\tau_{|M}.$ 
    \item $(\alpha_t)_{t\in\mathbb R}\subset {\text Aut}(\tilde M,\tilde\tau)$ is a $1$-parameter group with $\lim_{t\to 0}\|\alpha_t(x)-x\|_2=0$, for any $x\in\tilde M.$ 
    \item There is $\beta\in {\text Aut}(\tilde M,\tilde\tau)$ satisfying $\beta_{|M}=\text{Id}_M$, $\beta^2=\text{Id}_{\tilde M}$ and $\beta\alpha_t=\alpha_{-t}\beta$, for any $t\in\mathbb R$.
    
          \item $\alpha_t$ does not converge uniformly to the identity on $(M)_1$ as $t\to 0$.

\end{itemize}

\end{definition}

For a subalgebra $Q\subset qMq$, we say that $Q$ is {\it $\alpha$-rigid} if $\alpha_t$  converges uniformly to the identity on $(Q)_1$ as $t\to 0$.
We will repeatedly use the following stability result for s-malleable deformations.  

\begin{proposition}\emph{\cite[Proposition 3.4]{Va10b}}\label{proposition.deformation.set.vaes}
Let $(\tilde M, (\alpha_t)_{t\in\mathbb R})$ be an s-malleable deformation of a tracial von Neumann algebra $M$. Let $P\subset pMp$ be a subalgebra that is generated by a group of unitaries $\mathcal G\subset \mathcal U(P)$. Assume that $\alpha_t\to {\rm id}$ uniformly on $r\mathcal G r$ for a projection $r\in pMp$.

Then there is a projection $z\in \mathcal N_{pMp}(P)'\cap pMp $ with  $r\leq z$ such that
$Pz$ is $\alpha$-rigid. 
\end{proposition}

\subsection{Definition of class $\mathscr{M}$}
We are now ready to define the class of II$_1$ factors that is used in our main results stated in the introduction. 

\begin{definition}\label{definition.classM}
We say that a non-amenable II$_1$ factor $M$ belongs to class $\bm{\mathscr{M}}$ if there exists an s-malleable deformation $(\tilde M, (\alpha_t)_{t\in\mathbb R})$ of $M$ and an amenable subalgebra $A\subset M$ satisfying $L^2(\tilde M)\ominus L^2(M)$ is a mixing $M$-$M$-bimodule relative to $A$, 
         $L^2(\tilde M)\ominus L^2(M)$ is weakly contained in the coarse bimodule $L^2(M)\otimes L^2(M)$ as $M$-$M$-bimodules and one of    
    the following holds: 
    \begin{enumerate}
        \item $A=\mathbb C1$;
        
        \item If $N$ is a tracial von Neumann algebra and  $P\subset p(M\bar\otimes N)p$ a subalgebra such that $P\prec_{M\otimes N} A\otimes N$ and  $P'\cap p(M\bar\otimes N)p$ is strongly non-amenable relative to $1\otimes N$, then $P\prec_{M\otimes N} 1\otimes N.$

    
\end{enumerate}

As a consequence of Popa's spectral gap principle \cite{Po06b}, we continue with the following remark.

\begin{remark}\label{remark.popa.spectral.gap}
Let $M\in\bm{\mathscr{M}}$ be a II$_1$ factor, denote by $(\tilde M, (\alpha_t)_{t\in\mathbb R})$ the associated s-malleable deformation of $M$ and let $N$ be any tracial von Neumann algebra. If $Q\subset q(M\bar\otimes N)q$ is a von Neumann subalgebra which is strongly non-amenable relative to $N$, then $Q'\cap q(M\bar\otimes N)q$ is $(\alpha\otimes {\rm id})$-rigid.
This essentially follows from  \cite{Po06b} (see, for instance, the proofs of \cite[Lemma 2.2]{Io11} or \cite[Lemma 3.5]{Dr20}).
\end{remark}

\begin{proposition}\label{remark.classM.classA}
If $\Gamma$ belongs to one of the three classes of groups mentioned before Theorem \ref{A}, then $L(\Gamma)$ belongs to $\bm{\mathscr{M}}$. This follows, for instance, from Remark \ref{remark.popa.spectral.gap} and \cite[Proposition 3.4]{Dr20}).
\end{proposition}

Next, we present a useful result for group von Neumann algebras $L(\Gamma)$ that belong to $\bm{\mathscr M}$ in order to understand structural results of trace preserving actions of $\Gamma$. In fact, this is a direct consequence of Popa's spectral gap principle \cite{Po06b}.

\begin{lemma}\label{lemma.trick.commultiplication.psi}
Let $\Gamma\car B$ be a trace preserving action and denote  $\mathcal M=B\rtimes\Gamma$. 
We denote by $\Psi:\mathcal M\to \mathcal M \bar\otimes L(\Gamma)$  the $*$-homomorphism given by $\Psi(bu_g)=bu_g\otimes u_g$, for all $b\in B$ and $g\in\Gamma$.
Assume that $L(\Gamma)$ belongs to $\bm{\mathscr{M}}$ and let $(\tilde M, (\alpha_t)_{t\in\mathbb R})$ be the associated s-malleable deformation. 

If $P\subset p \mathcal M p$ be a von Neumann subalgebra that is strongly non-amenable relative to $B$ inside $
\mathcal M$, then $\Psi (P'\cap p\mathcal M p)$ is ${\rm (id}\otimes \alpha)$-rigid. Moreover,  if we also assume that $P'\cap p\mathcal M p\nprec_{\mathcal M} B$, then $\Psi (P\vee (P'\cap p\mathcal M p))$ is ${\rm (id}\otimes \alpha)$-rigid.
\end{lemma} 
 
{\it Proof.} 
Denote $M=L(\Gamma)$.
By applying \cite[Lemma 2.10]{Dr19a}, we get that $\Psi(P)$ is strongly non-amenable relative to $\mathcal M\bar\otimes 1$ inside $\mathcal M\bar\otimes M$. Then Remark \ref{remark.popa.spectral.gap} implies that $\Psi (P'\cap p\mathcal M p)$ is ${\rm (id}\otimes \alpha)$-rigid. 
For proving the second part, assume in addition that $P'\cap p\mathcal M p\nprec_{\mathcal M} B$ and let $A\subset M$ be an amenable subalgebra as given by the assumption that $M$ belongs to $\bm{\mathscr M}$. By using  \cite[Lemma 9.2(1)]{Io11} we get $\Psi(P'\cap p\mathcal M p)\nprec_{\mathcal M\bar\otimes M} \mathcal M\otimes 1$, and by assumption we must have $\Psi(P'\cap p\mathcal M p)\nprec_{\mathcal M\bar\otimes M} \mathcal M\otimes A$. 
Hence, in combination with $\Psi (P'\cap p\mathcal M p)$ is ${\rm (id}\otimes \alpha)$-rigid
 and the fact that $L^2(\mathcal M\bar\otimes \tilde M)\ominus L^2(\mathcal M\bar\otimes M)$ is a mixing $\mathcal M\bar\otimes M$-bimodule relative to $\mathcal M\bar\otimes A$,
 we get from \cite[Corollary 6.7]{dSHHS20} that $\Psi (P\vee (P'\cap p\mathcal M p))$ is ${\rm (id}\otimes \alpha)$-rigid.
\hfill$\blacksquare$

\end{definition}

\subsection{Measure equivalence and non property Gamma for class $\mathscr{M}$}

In this subsection we show that the  lack of property Gamma is preserved under measure equivalence for finite products of groups whose von Neumann algebras belong to $\bm{\mathscr M}$, see Proposition \ref{property.gamma}. For proving this result, we first establish the following notation that will be assumed for Proposition \ref{property.gamma}, but
will also be useful for the following sections.

\begin{notation}\label{notation} Let $\Lambda$ be a countable icc group that is measure equivalent to a product $\Gamma=\Gamma_1\times\dots\times\Gamma_n$ of $n\ge 1$ groups. By using \cite[Lemma 3.2]{Fu99}, there exist $d\ge 1$, free ergodic pmp actions $\Gamma\car (X,\mu)$ and $\Lambda\car (Y,\nu)$ such that 
$$
\mathcal R(\Lambda\car Y)=\mathcal R(\Gamma\times \mathbb Z/d\mathbb Z \car X \times \mathbb Z/d\mathbb Z)\cap (Y\times Y).
$$
Here, we considered that $\mathbb Z/d\mathbb Z \car (\mathbb Z/d\mathbb Z,c)$ acts by addition and $c$ is the counting measure. We also identified $Y$ as a measurable subset of $X\times \mathbb Z/d\mathbb Z$ and denote $p=1_Y \in L^\infty(X\times \mathbb Z/d\mathbb Z)$. Note that $L^\infty(\mathbb Z/d\mathbb Z)\rtimes \mathbb Z/d\mathbb Z=\mathbb M_d(\mathbb C)$. Hence, by letting $B=L^\infty(Y)$, $A=L^\infty(X)\otimes M_d(\mathbb C)$ and $M=A\rtimes\Gamma$, we have $pMp=B\rtimes\Lambda$ and $B\subset pAp$. Denote by $\{u_g\}_{g\in \Gamma}$ and $\{v_\lambda\}_{\lambda\in \Lambda}$ the canonical unitaries implementing the actions $\Gamma\car A$ and $\Lambda\car B$, respectively.

Following \cite{PV09} we define the $*$-homomorphism $\Delta: pMp\to pMp \bar\otimes L(\Lambda)$ by $\Delta(bv_\lambda)=bv_\lambda\otimes v_\lambda$, for all $b\in B,\lambda\in\Lambda$. One can extend $\Delta$ to a $*$-homomorphism $\Delta:M\to M\bar\otimes L(\Lambda)$ and verify that  $\Delta (M)'\cap M\bar\otimes L(\Lambda)=\mathbb C 1$ since $\Lambda$ is icc (see the first part of \cite[Section 5]{DHI16} for more details). 

For any $i\in\overline{1,n}$, let $\Psi^i: M\to M\bar\otimes L(\Gamma_i)$ be the $*$-homomorphism given by $\Psi^i(xu_g)=xu_g\otimes u_g$, for all $x\in A\rtimes\Gamma_{\widehat i} ,g\in\Gamma_i$.

\end{notation}






\begin{proposition}\label{property.gamma}
Assume that $L(\Gamma_i)$ belongs to $\bm{\mathscr M}$, for any $1\leq i\leq n$.

Then $L(\Lambda)$ does not have property Gamma.
\end{proposition}

{\it Proof.} Let $(\tilde M, (\alpha^i_t)_{t\in\mathbb R})$ be the associated s-malleable deformation of $L(\Gamma_i) \in\bm{\mathscr{M}}$. By assuming that $L(\Lambda)$ has property Gamma, we can use \cite[Theorem 3.1]{HU15} to obtain a decreasing sequence of diffuse abelian von Neumann subalgebras $Q_n\subset L(\Lambda)$ with $n\ge 1$ such that $L(\Lambda)=\bigvee_{n\ge 1} (Q_n'\cap L(\Lambda)$). Since $Q_1$ is abelian note that for all $k\leq n$, we have
\begin{equation}\label{dd}
    \mathcal Z(Q_n'\cap pMp)\subset Q_1'\cap pMp\subset Q_k'\cap pMp,
\end{equation}

Using Zorn's lemma and a maximality argument, one can show that for any $m\ge 1$, there exist maximal projections $r_m^1,\dots, r_m^n\in Q_m'\cap pMp$ satisfying $Q_mr_m^i\nprec_{M}A\rtimes \Gamma_{\widehat i}$, for any $i\in\overline{1,n}$. One can check that $r_m^i\in \mathcal Z(Q_m'\cap pMp)$ and $Q_m (p-r_m^i)\prec_M^s A\rtimes\Gamma_{\widehat i} $, for any $i\in\overline{1,n}$ (see the proof of \cite[Lemma 2.4]{DHI16}).

Since $Q_{m}\nprec_{M} A$, \cite[Lemma 2.8(2)]{DHI16} implies that $\bigwedge _{i=1}^n (p-r_m^i)=0$, which proves that $\bigvee _{i=1}^n r_m^i =p$.
Hence, for any $m\ge 1$ there is $i_m\in \overline{1,n}$ such that $\tau(r_m^{i_m})\ge \frac{\tau(p)}{n}$. Up to passing to a subsequence, we can assume that there is $j\in \overline{1,n}$ such that $i_m=j$, for all $m\ge 1$.
Next, note that \eqref{dd} gives that $r_{m}^j\in \mathcal Z(Q_m'\cap pMp)\subset Q_{m-1}'\cap pMp$. Since $Q_{m-1} r_m^j \nprec_M A\rtimes \Gamma_{\widehat j}$, it follows from the choice of all the $r_m^j$'s that  $\{r_{m}^j\}_{m\ge 1}$ is a decreasing sequence of projections. If we let $r^j=\bigwedge_{m\ge 1} r_{m}^j$, we deduce that $r^j$ is a non-zero projection since $\tau(r^j)\ge \frac{\tau(p)}{n}$. For all $m\ge k\ge 1$, since $Q_m\subset Q_k$ we have that $r_m^j \in (Q_k'\cap pMp)'\cap pMp$. Consequently, by letting $m\to\infty$, we deduce that $r^j \in (Q_k'\cap pMp)'\cap pMp$, for all $k\ge 1$, which implies  that $r^j\in L(\Lambda)'\cap pMp=\mathbb C p$. Since $r^j\neq 0$, we derive that $r^j=p$, and therefore, we must have $r_m^j=p$, for any $m\ge 1$. This implies that $Q_m\nprec_M A\rtimes\Gamma_{\widehat j}$, for any $m\ge 1$.

Since $\Gamma_j$ is non-amenable, it follows from \cite[Lemma 2.9]{DHI16} that $L(\Lambda)$ is non-amenable relative to $A\rtimes \Gamma_{\widehat j}$ inside $M$.
Since relative amenability is closed under inductive limits (see \cite[Lemma 2.7]{DHI16}), there exists $k\ge 1$ such that $Q_k'\cap pMp$  is non-amenable relative to  $A\rtimes \Gamma_{\widehat j}$ inside  $M$.    
Using \cite[Lemma 2.6]{DHI16} there is a non-zero projection $z^j\in\mathcal Z(Q_k'\cap pMp)$ such that  
\begin{equation}\label{dd1}
(Q_k'\cap pMp)z^j \text{ is strongly non-amenable relative to } A\rtimes \Gamma_{\widehat j} \text{ inside  }M.    
\end{equation}
This implies by Lemma \ref{lemma.trick.commultiplication.psi} that $\Psi^j (Q_k z^j) \text{ is }{\rm (id}\otimes \alpha^j) \text{-rigid}$. Fix an arbitrary $m\ge k$. Since $Q_m\subset Q_k$,  we have $z^j\in\mathcal Z(Q_k'\cap pMp)\subset Q_m'\cap pMp$ and  
\begin{equation}\label{dd2}
\Psi^j (Q_m z^j) \text{ is }{\rm (id}\otimes \alpha^j) \text{-rigid} \text{  and  }Q_mz^j\nprec_M A\rtimes\Gamma_{\widehat j}.
\end{equation}
Equation \eqref{dd1} also implies that
\begin{equation}\label{dd3}
z_j(Q_m'\cap pMp)z^j \text{ is strongly non-amenable relative to } A\rtimes \Gamma_{\widehat j} \text{ inside  }M.    
\end{equation}
By combining \eqref{dd2} and \eqref{dd3}, it follows from the second part of Lemma \ref{lemma.trick.commultiplication.psi} that $\Psi^j (z^j(Q_m'\cap pMp) z^j) \text{ is }{\rm (id}\otimes \alpha^j) \text{-rigid}$, for any $m\ge k$. Note that \eqref{dd2} gives in particular that $z^j(Q_m'\cap pMp) z^j\nprec_M A\rtimes\Gamma_{\widehat j}$, for any $m\ge k$. Therefore, we may 
apply \cite[Theorem 3.5]{dSHHS20} (see also \cite[Theorem 3.2]{Dr20}) to deduce that $\Psi^j (z^ j \bigvee_{m\ge k}(Q_m'\cap pMp) z^j) \text{ is }{\rm (id}\otimes \alpha^j) \text{-rigid}$.
Using \cite[Proposition 5.6]{dSHHS20} there exists a non-zero projection $\tilde z^j\in  \mathcal Z(\bigvee_{m\ge n}(Q_m'\cap pMp))$ such that  $\Psi^j ( \bigvee_{m\ge k}(Q_m'\cap pMp) \tilde z^j) \text{ is }{\rm (id}\otimes \alpha^j) \text{-rigid}$. Note however that $\bigvee_{m\ge k}(Q_m'\cap pMp)$ is a factor since $\bigvee_{m\ge k}(Q_k'\cap L(\Lambda))=L(\Lambda)$ and $\Lambda$ is icc.
In particular, $\Psi^j (L(\Lambda) )$ is ${(\rm id}\otimes \alpha^j)$-rigid. Since $\Psi^j(B)\subset M\otimes 1$, it follows that $\Psi^j (M) )$ is ${(\rm id}\otimes \alpha^j)$-rigid, which gives that $L(\Gamma_j)$ is $\alpha^j$-rigid, contradiction. This ends the proof of the proposition.
\hfill$\blacksquare$

\section{Measure equivalence and tensor product decompositions for class $\mathscr{M}$}

In this section we establish the main ingredients needed for the proof of Theorem \ref{A} by building upon methods from \cite{DHI16,IM19}. Throughout this section we will use Notation \ref{notation} and the following assumption.

\begin{assumption}\label{assumption}
For any $i\in\overline{1,n}$, assume that $L(\Gamma_i)$ belongs to $\bm{\mathscr M}$ and denote by $(\tilde M, (\alpha^i_t)_{t\in\mathbb R})$  the associated s-malleable deformation of $L(\Gamma_i)$. 
\end{assumption}

\subsection{Step 1}
The main goal of this subsection is to prove Theorem \ref{Step1}. For doing this, we will use a relative version of the flip automorphism method introduced by Isono and Marrakchi \cite{IM19}.

\begin{theorem}\label{Step1}
Let $L(\Lambda)=P_1\bar\otimes P_2$ be a tensor product decomposition into II$_1$ factors.

Then there is a partition $S_1\sqcup S_2=\{1,\dots,n\}$ into non-empty sets such that $\Delta(A\rtimes\Gamma_{S_i})\prec^s_{M\bar\otimes L(\Lambda)} M\bar\otimes P_i$, for all $i\in \{1,2\}$.
\end{theorem}

Before proceeding to the proof of Theorem \ref{Step1}, we need the following two lemmas.

\begin{lemma}\label{lemma.first.step}
Let $L(\Lambda)=P_1\bar\otimes P_2$ be a tensor product decomposition into II$_1$ factors.

Then there is a partition $S_1\sqcup S_2=\{1,\dots,n\}$ and a projection $0\neq z\in \Delta(L(\Gamma))'\cap \mathcal M$ such that:
\begin{itemize}
    \item $\Delta(L(\Gamma_i))z$ is strongly non-amenable relative to $M\bar\otimes P_{1}$ inside $\mathcal M$ for all $i\in S_2.$
    \item $\Delta(L(\Gamma_i))z$ is strongly non-amenable relative to $M\bar\otimes P_{2}$ inside $\mathcal M$ for all $i\in S_1.$
\end{itemize}
\end{lemma}

{\it Proof.} Let  $i\in\{1,\dots,n\}$. Since $\Gamma_i$ is non-amenable, by \cite[Proposition 2.4]{KV15} we get that $\Delta(L(\Gamma_i))$ is strongly non-amenable relative to $M\otimes 1$ inside $\mathcal M$. It follows that for every non-zero projection  $z\in \Delta(L(\Gamma))'\cap \mathcal M$, there exist $f(i,z)\in\{1,2\}$ and a non-zero projection $p(i,z)\in \Delta(L(\Gamma))'\cap \mathcal M$ with $p(i,z)\leq z$ such that
\begin{equation}\label{ddd}
\Delta(L(\Gamma_i))p(i,z) \text{ is strongly non-amenable relative to } M\bar\otimes P_{f(i,z)} \text{  inside } \mathcal M.     
\end{equation}
Indeed, otherwise there exists a non-zero projection $z\in \Delta(L(\Gamma))'\cap \mathcal M$ such that for any $k\in\{1,2\}$ and non-zero projection $z_0\in \Delta(L(\Gamma))'\cap \mathcal M$ with $z_0\leq z$, there exists a non-zero projection $\tilde z_0\in \Delta(L(\Gamma))'\cap \mathcal M$ with $\tilde z_0\leq z_0$ for which
$\Delta(L(\Gamma_i))\tilde z_0$ is amenable relative to $M\bar\otimes P_k$ inside $\mathcal M$.
By using \cite[Proposition 2.7]{PV11} we derive that there exists a non-zero projection $\tilde z_1\in \Delta(L(\Gamma))'\cap \mathcal M$ with $\tilde z_1\leq z$ for which
$\Delta(L(\Gamma_i))\tilde z_1$ is amenable relative to $M\bar\otimes 1$ inside $\mathcal M$, contradiction.

By applying \eqref{ddd} finitely many times, the proof will be obtained as follows.
Define $z_1=p(1,1)$ and $f(1)=f(1,1)$. For any $i\in\{2,\dots,n\}$ we recursively define $z_i=p(i,z_{i-1})$ and $f(i)=f(i,z_{i-1})$. Note that $z_1\ge z_2\ge\dots\ge z_n$ are non-zero projections in $\Delta(L(\Gamma))'\cap \mathcal M$. Hence, the lemma follows by letting  $S_1=f^{-1}(2), S_2=f^{-1}(1)$ and $z=z_n$. 
\hfill$\blacksquare$

We continue with the following notation that will be used in  the following lemma, but also in the proof of Theorem \ref{Step1}.  For any $1\leq j\leq n$, denote $\Psi^{j,4}={\rm id}\otimes {\rm id}\otimes {\rm id}\otimes \Psi^j$ and $\alpha^{j,5}= {\rm id}\otimes{\rm id}\otimes {\rm id}\otimes {\rm id}\otimes \alpha^j$. By letting $\mathcal M=M\bar\otimes L(\Lambda)$, note  that $\Psi^j(p)=p\otimes 1$ and
$
\Psi^{j,4}(\mathcal M\bar\otimes \mathcal M)\subset M\bar\otimes L(\Lambda) \bar\otimes M \bar\otimes pMp \bar\otimes L(\Gamma_j).
$

\begin{lemma}\label{lemma.deformation.not.convergence}
Let $\sigma\in {\rm Aut}(\mathcal M\bar\otimes \mathcal M)$ be an automorphism for which $\sigma_{|(M\otimes 1)\bar\otimes (M\otimes 1))}={\rm Id}$ and $(1\otimes L(\Lambda))\bar\otimes (1\otimes L(\Lambda))$ is $\sigma$-invariant. Then $\Psi^{j,4}(\sigma(\Delta(M)\bar\otimes \Delta(M))z$ is not $\alpha^{j,5}$-rigid, for all non-zero projections $z\in \Psi^{j,4}(\sigma(\Delta(M)\bar\otimes \Delta(M))'\cap M\bar\otimes L(\Lambda) \bar\otimes M \bar\otimes pMp \bar\otimes L(\Gamma_j)$ and $j\in\overline{1,n}$.
\end{lemma}

{\it Proof.} 
By assuming the contrary, there exist  $j\in\{1,\dots,n\}$ and a projection $z$ as in the statement such that $\Psi^{j,4}(\sigma(\Delta(M)\bar\otimes \Delta(M))z$ is  $\alpha^{j,5}$-rigid.
Hence, for any $\epsilon>0$, there is $t_0>0$ such that
$$
\| \Psi^{j,4}(\sigma( v_g\otimes v_g\otimes v_h\otimes v_h ))z - \alpha_t^{j,5}( \Psi^{j,4}(\sigma( v_g\otimes v_g\otimes v_h\otimes v_h ))z)  \|_2\leq \epsilon,
$$
for all $g,h\in\Lambda$ and $|t|\leq t_0.$ Since $\sigma$ acts trivially on $( M\otimes 1)\bar\otimes ( M \otimes 1)$, we obtain that
$$
\| \Psi^{j,4}(\sigma( 1\otimes v_g\otimes 1\otimes v_h ))z - \alpha_t^{j,5}( \Psi^{j,4}(\sigma( 1\otimes v_g\otimes 1\otimes v_h ))z)  \|_2\leq \epsilon,
$$
for all $g,h\in\Lambda$ and $|t|\leq t_0.$ If we let $\mathcal G=\{\Psi^{j,4}(\sigma( 1\otimes v_g\otimes 1\otimes v_h ))|\; g,h\in\Lambda\}$, we get that $\mathcal G''=\Psi^{j,4}(1\otimes L(\Lambda)\otimes 1\bar\otimes L(\Lambda))=1\otimes L(\Lambda)\otimes 1\bar\otimes \Psi^j(L(\Lambda)) $ since $(1\otimes L(\Lambda))\bar\otimes (1\otimes L(\Lambda))$ is $\sigma$-invariant. Note that $\mathcal N_{M\bar\otimes L(\Lambda)\bar\otimes M\bar\otimes pMp\bar\otimes L(\Gamma_j)}(\mathcal G'')\subset 1\otimes 1\otimes 1\otimes (\Psi^{j}(L(\Lambda))'\cap (pMp\bar\otimes L(\Gamma_j)))$.

By applying Proposition \ref{proposition.deformation.set.vaes} we obtain a non-zero projection $z_0\in \Psi^{j}(L(\Lambda))'\cap (pMp\bar\otimes L(\Gamma_j))$ such that $\Psi^{j}(L(\Lambda))z_0$ is $({\rm id}\otimes\alpha^j)$-rigid. Since $\Psi^j(B)=B\otimes 1$ and $\Psi^{j}(pMp)'\cap (pMp\bar\otimes L(\Gamma_j))=\mathbb C (p\otimes 1)$, it follows from Proposition \ref{proposition.deformation.set.vaes} that $\Psi^{j}(pMp)$ is $({\rm id}\otimes\alpha^j)$-rigid, and hence, $\Psi^{j}(M)$ is $({\rm id}\otimes\alpha^j)$-rigid. This shows that  $L(\Gamma_j)$ is $\alpha^j$-rigid, contradiction.
\hfill$\blacksquare$

{\it Proof of Theorem \ref{Step1}.} 
 Denote $\mathcal M=M\bar\otimes L(\Lambda)$ and $\tilde{\mathcal M}=M\bar\otimes M$. For proving this theorem, we use the following variation of the flip automorphism method  from \cite{IM19}. Namely, since $L(\Lambda)=P_1\bar\otimes P_2$, we define $\sigma_{P_1}\in {\rm Aut}(\mathcal M\bar\otimes \mathcal M)$ by letting $\sigma_{P_1}(m\otimes p_1\otimes p_2 \otimes m'\otimes p_1'\otimes p_2')=m\otimes p_1'\otimes p_2 \otimes m'\otimes p_1\otimes p_2'$, for all $m,m'\in M,p_1,p_1'\in P_1, p_2,p_2'\in P_2$.


By applying Lemma \ref{lemma.first.step} and Lemma \ref{lemma.tensor} we obtain a partition $S_1\sqcup S_2=\{1,\dots,n\}$ and a non-zero projection $z\in \Delta(L(\Gamma))'\cap \mathcal M$ such that 
\begin{equation}\label{s1}
    \begin{array}{cc}
         &  \Delta(L(\Gamma_i))z\otimes 1 \text{ is strongly non-amenable relative to }(M\bar\otimes P_2)\bar\otimes (M\bar\otimes P_1) \text{ inside } \mathcal M\bar\otimes\mathcal M, \\
         &  1\otimes \Delta(L(\Gamma_j))z \text{ is strongly non-amenable relative to }(M\bar\otimes P_2)\bar\otimes (M\bar\otimes P_1) \text{ inside } \mathcal M\bar\otimes\mathcal M,
    \end{array}
\end{equation}
for all $i\in S_1$ and $j\in S_2$. By applying the flip automorphism $\sigma_{P_1}$ to \eqref{s1}, we derive that
\begin{equation*}\label{s2}
    \begin{array}{cc}
         &  \sigma_{P_1}(\Delta(L(\Gamma_i))z\otimes 1) \text{ is strongly non-amenable relative to }\mathcal M\bar\otimes (M\otimes 1) \text{ inside } \mathcal M\bar\otimes\mathcal M, \\
         &  \sigma_{P_1}(1\otimes \Delta(L(\Gamma_j))z) \text{ is strongly non-amenable relative to }\mathcal M\bar\otimes (M\otimes 1) \text{ inside } \mathcal M\bar\otimes\mathcal M,
    \end{array}
\end{equation*}
for all $i\in S_1$ and $j\in S_2$. By using Lemma \ref{lemma.tensor}, we further deduce  that 
\begin{equation}\label{s3}
    \begin{array}{cc}
         &  \sigma_{P_1}(\Delta(L(\Gamma_i))z\otimes 1) \text{ is strongly non-amenable relative to }\mathcal M\bar\otimes (M\otimes 1) \text{ inside } \mathcal M\bar\otimes\tilde{\mathcal M}, \\
         &  \sigma_{P_1}(1\otimes \Delta(L(\Gamma_j))z) \text{ is strongly non-amenable relative to }\mathcal M\bar\otimes (M\otimes 1) \text{ inside } \mathcal M\bar\otimes\tilde{\mathcal M},
    \end{array}
\end{equation}
for all $i\in S_1$ and $j\in S_2$. Denote $\widehat z=\sigma_{P_1}(z\otimes z)\in \sigma_{P_1}(\Delta(L(\Gamma))\bar\otimes \Delta(L(\Gamma)))'\cap (\mathcal M\bar\otimes \mathcal M)$ and note that $\widehat z\leq \sigma_{P_1}(z\otimes 1)$, $\widehat z\leq \sigma_{P_1}(1\otimes z)$. For ease of notation, we denote $Q_i=\Delta(L(\Gamma_i))\otimes 1$ and $R_i=\Delta(L(\Gamma_{\widehat i})\bar\otimes \Delta(L(\Gamma))$, for all $i\in S_1$. Similarly, denote $Q_j=1\otimes \Delta(L(\Gamma_j))$ and $R_j=\Delta(L(\Gamma))\bar\otimes \Delta(L(\Gamma_{\widehat j}) $, for all $j\in S_2$. Note that $Q_i \vee R_i=\Delta(L(\Gamma))\bar\otimes \Delta(L(\Gamma))$, for any $i\in\{1,\dots,n\}$.

By applying a similar argument to the one used in Lemma \ref{lemma.first.step} we deduce from \eqref{s3} that there exist a non-zero projection $\tilde z\in \sigma_{P_1}(\Delta(L(\Gamma))\bar\otimes \Delta(L(\Gamma)))'\cap (\mathcal M\bar\otimes \tilde{\mathcal M})$ and a function $\varphi:\{1,\dots,n\}\to \{1,\dots,n\}$ such that
\begin{equation}\label{s4}
    \sigma_{P_1}(Q_i)\tilde z \text{ is strongly non-amenable relative to } \mathcal M\bar\otimes M\bar\otimes (A\rtimes \Gamma_{\widehat{\varphi(i)}}),
\end{equation}
for all $i\in\{1,\dots,n\}$. By Lemma \ref{lemma.trick.commultiplication.psi}, we get that for any $i\in\{1,\dots,n\}$,
\begin{equation}\label{s5}
    \Psi^{g(i),4}(\sigma_{P_1}(R_i)\tilde{z})  \text{ is }  \alpha^{g(i),5} \text{ -rigid}.
\end{equation}

Next, we claim that the map $\varphi$ is bijective. If this does not hold, it is easy to see that we can deduce from \eqref{s4} that there exists $j\in\{1,\dots,n\}$ such that 
\begin{equation}\label{ss5}
\Psi^{j,4}(\sigma_{P_1}(\Delta(L(\Gamma))\bar\otimes \Delta(L(\Gamma)))\tilde{z}) \text{ is  }\alpha^{j,5}\text{-rigid}.    
\end{equation}

Using the position of $B\subset A$, and $\sigma_{P_1}(\Delta(B)\bar\otimes \Delta(B))\subset \mathcal M\bar\otimes (M\otimes 1)$, we obtain that
\begin{equation}\label{sss5}
\Psi^{j,4}(\sigma_{P_1}(\Delta(A)\bar\otimes \Delta(A)) \text{  is } \alpha^{j,5} \text{-rigid}.     
\end{equation}

Relations \eqref{ss5} and \eqref{sss5} in combination with
Proposition \ref{proposition.deformation.set.vaes} gives a contradiction to Lemma \ref{lemma.deformation.not.convergence}. This shows that $\varphi$ is indeed bijective. Next, we claim that for all $i\in\{1,\dots,n\}$,
\begin{equation}\label{s6}
    \sigma_{P_1}(R_i)\tilde{z}\prec^s_{\mathcal M\bar\otimes \tilde{\mathcal M}} \mathcal M\bar\otimes M\bar\otimes (A\rtimes \Gamma_{\widehat{g(i)}}),
\end{equation}
Assume by contradiction that there is $i\in\{1,\dots, n\}$ for which \eqref{s6} does not hold. Then by using \cite[Lemma 2.4(2)]{DHI16}, it follows that, up to replacing $\tilde z$ by a smaller non-zero projection, we have $\sigma_{P_1}(R_i)\tilde{z}\nprec_{\mathcal M\bar\otimes \tilde{\mathcal M}} \mathcal M\bar\otimes M\bar\otimes (A\rtimes \Gamma_{\widehat{g(i)}}).$
Using \eqref{s4} and \eqref{s5} we may apply Lemma \ref{lemma.trick.commultiplication.psi} to deduce that
$\Psi^{j,4}(\sigma_{P_1}(\Delta(L(\Gamma))\bar\otimes \Delta(L(\Gamma)))\tilde{z}) \text{ is  }\alpha^{j,5}\text{-rigid}.$ As before, Proposition \ref{proposition.deformation.set.vaes} leads to a contradiction.

Finally, by applying \cite[Lemma 2.8(2)]{DHI16} finitely many times, we deduce from \eqref{s6} that
$$\sigma_{P_1}(\Delta(L(\Gamma_{\widehat{S_1}}))\bar\otimes \Delta(L(\Gamma_{\widehat{S_2}}))) \prec_{\mathcal M\bar\otimes \tilde{\mathcal M}} \mathcal M\bar\otimes (M \otimes 1).$$ 
By applying Lemma \ref{lemma.tensor} we further obtain that $\sigma_{P_1}(\Delta(L(\Gamma_{\widehat{S_1}}))\bar\otimes \Delta(L(\Gamma_{\widehat{S_2}}))) \prec_{\mathcal M\bar\otimes {\mathcal M}} \mathcal M\bar\otimes (M \otimes 1)$. By applying the flip automorphism $\sigma_{P_1}$ to the previous intertwining relation,
we deduce that $\Delta(L(\Gamma_{S_i}))\prec_{\mathcal M} M\bar\otimes P_i$, for all $i\in \{1,2\}$. Since $\Delta(A)\prec^s_{\mathcal M}B\otimes 1$, we may use \cite[Lemma 2.3]{BV12} to get that $\Delta(A\rtimes\Gamma_{S_i})\prec_{\mathcal M} M\bar\otimes P_i$, for all $i\in \{1,2\}$. Since $\mathcal N_{\mathcal M}(\Delta(A\rtimes\Gamma_{S_i}))'\cap \mathcal M\subset \Delta(M)'\cap\mathcal M=\mathbb C1$, we obtain $\Delta(A\rtimes\Gamma_{S_i})\prec^s_{\mathcal M} M\bar\otimes P_i$, for all $i\in \{1,2\}$.

For showing that $S_1$ and $S_2$ are non-empty sets, we suppose the contrary. Hence, without loss of generality assume that $S_2$ is empty. This shows that $\Delta(M)\prec_{\mathcal M} M\bar\otimes P_1$, which implies from \cite[Lemma 9.2]{Io10} that $L(\Lambda)\prec_{L(\Lambda)} P_1$. This shows that $P_2$ is not diffuse, contradiction.
\hfill$\blacksquare$

\subsection{Step 2.} By using Step 1, we obtain the following intertwining result. Recall that we are using Notation \ref{notation} and Assumption \ref{assumption}.

\begin{theorem}\label{Step2} Let $L(\Lambda)=P_1\bar\otimes P_2$ be a tensor product decomposition into II$_1$ factors.

Then there is a partition $T_1\sqcup T_2=\{1,\dots,n\}$ such that $P_i\prec^s_{M} A\rtimes \Gamma_{T_i}$, for all $i\in\{1,2\}.$
\end{theorem}

Throughout the proof we are using the following notation: if $N$ is a tracial von Neumann algebra and $P\subset pNp$ and $Q\subset qNq$ are von Neumann subalgebras, we denote $P\prec_N^{s'}Q$ if $P\prec_N Qq'$, for any non-zero projection $q'\in Q'\cap qNq$.

{\it Proof of Theorem \ref{Step2}.} 
Theorem \ref{Step1} implies that there exist projections $r_1\in \Delta(L(\Gamma_{S_1})), q_1\in M\bar\otimes P_1$, a non-zero partial isometry $w_1\in q_1 (M\bar\otimes M) r_1$ and a $*$-homomorphism $\varphi_1: r_1\Delta(L(\Gamma_{S_1}))r_1\to q_1(M\bar\otimes P_1)q_1$ such that $\varphi_1(x)w_1=w_1x$, for all $x\in r_1\Delta(L(\Gamma_{S_1}))r_1$. Fix an arbitrary $j_0\in S_1$. Since $L(\Gamma_{j_0})$ is a II$_1$ factor we can  apply \cite[Lemma 4.5]{CdSS17} and therefore assume without loss of generality that $r_1\in L(\Gamma_{j_0})$. In addition, we can assume that the support projection of $E_{M\bar\otimes P_1}(w_1w_1^*)$ equals $q_1$. For any $j\in S_1$,
denote $Q_1^j=\varphi_1 (r_1\Delta(L(\Gamma_j))r_1)\subset q_1(M\bar\otimes P_1)q_1$ and let $Q_1=\bigvee_{j\in S_1} Q_1^j$. Note that for any subset $S\subset S_1$, we have $\Delta(L(\Gamma_S))\prec^{s'}_{M\bar\otimes P_1} \bigvee_{j\in S}Q_1^j$. Indeed, let $S\subset S_1$ and consider a non-zero projection $z\in Q_1'\cap q_1 (M\bar\otimes P_1) q_1$. Note that $\tilde w_1:=zw_1\neq 0$ since otherwise $zE_{M\bar\otimes P_1}(w_1w_1^*)=0$, which implies that $z=0$, false. This shows that the $*$-homomorphism $\tilde\varphi_1:r_1\Delta(L(\Gamma_{S}))r_1\to \bigvee_{j\in S} Q_1^jz$ satisfies $\tilde\varphi_1 (x)\tilde w_1=\tilde w_1 x$, for all $x\in r_1\Delta(L(\Gamma_{S}))r_1.$ By replacing $\tilde w_1$ by the partial isometry from its polar decomposition, we derive that $\Delta(L(\Gamma_S))\prec_{M\bar\otimes P_1} \bigvee_{j\in S}Q_1^jz$. By using \cite[Lemma 2.4]{DHI16} it follows that  $\Delta(L(\Gamma_S))\prec^{s'}_{M\bar\otimes P_1} \bigvee_{j\in S}Q_1^j$. By applying \cite[Lemma 2.3]{Dr19b}, we derive that for any subset $S\subset S_1$,
\begin{equation}\label{y0}
\Delta(L(\Gamma_S))\prec^{s'}_{M\bar\otimes M} \bigvee_{j\in S} Q_1^j.    
\end{equation}
The rest of the proof is divided between three claims.

{\bf Claim 1.} For any $j\in S_1$ and non-zero projection $z\in Q_1'\cap q_1 (M\bar\otimes P_1) q_1$, there exist $k\in \{1,\dots,n\}$ and a non-zero projection $z_0 \in Q_1'\cap q_1 (M\bar\otimes P_1) q_1$ with $z_0\leq z$ such that
 $Q_1^{j} z_0$ { is strongly non-amenable relative to } $M\bar\otimes (A\rtimes\Gamma_{\widehat {k}})$ inside $M\bar\otimes M$.

{\it Proof of Claim 1.}
We assume by contradiction that there exist $j\in S_1$ and a non-zero projection $z\in Q_1'\cap q_1 (M\bar\otimes P_1) q_1$ such that $Q_1^{j} z$ {is amenable relative to} $M\bar\otimes (A\rtimes\Gamma_{\widehat {k}}),$ for all $k\in\{1,\dots,n\}$.
By applying \cite[Proposition 2.7]{PV11} we get that $Q_1^{j}z$ is amenable relative to $M\otimes 1$ inside $M\bar\otimes M$. 
Relation \eqref{y0} implies that $\Delta(L(\Gamma_j))\prec^{}_{M\bar\otimes M} Q_1^j z.$
We can apply \cite[Lemma 2.4(3) and Lemma 2.6(3)]{DHI16} and derive that there exists a non-zero projection $r'\in \Delta(L(\Gamma_{j}))'\cap M\bar\otimes M$ such that $\Delta(L(\Gamma_{j}))r'$ is amenable relative to $Q_1^{j}z\oplus \mathbb C(1-z)$. Using \cite[Proposition 2.4(3)]{OP07} we derive that $\Delta(L(\Gamma_{j}))r'$ is amenable relative to $M\otimes 1$. By using \cite[Lemma 10.2(5)]{IPV10} we deduce that $\Gamma_{j}$ is amenable, contradiction. Thus, there exist $k\in \{1,\dots, n\}$ such that $Q_1^{j} z$  is  non-amenable relative to  $M\bar\otimes (A\rtimes\Gamma_{\widehat {k}}).$ By \cite[Lemma 2.6]{DHI16}, there exists a non-zero projection $z_0 \in \mathcal N_{q_1(M\bar\otimes M)q_1}(Q_1^j)'\cap q_1 (M\bar\otimes M_1) q_1\subset Q_1'\cap   q_1(M\bar\otimes P_1)q_1$ with $z_0\leq z$ such that
 $Q_1^{j} z_0$ { is strongly non-amenable relative to } $M\bar\otimes (A\rtimes\Gamma_{\widehat {k}}).$
\hfill$\square$

By applying Claim 1 finitely many times and proceeding as in  the proof of Lemma \ref{lemma.first.step}, there exist a non-zero projection $z\in Q_1'\cap q_1 (M\bar\otimes P_1) q_1$ and a map $\overline{1,n}\ni j\to k_j\in\overline{1,n}$ such that
\begin{equation}\label{y1}
    Q_1^{j} z \text{ is strongly non-amenable relative to } M\bar\otimes (A\rtimes\Gamma_{\widehat {k_{j}}}), \text{ for any }j\in\overline{1,n}. 
\end{equation}

{\bf Claim 2.} $P_2\prec_M^s A\rtimes \Gamma_{\widehat {k_j}}$, for all $j\in S_1$.

{\it Proof of Claim 2.} Fix an arbitrary $j\in S_1$. We are in one of the following situations. First, by assuming that $(1\otimes P_2)z\prec_{M\bar\otimes M}M\bar\otimes (A\rtimes\Gamma_{\widehat{k_j}})$, we get $P_2\prec_M A\rtimes \Gamma_{\widehat {k_j}}$. Since $\mathcal N_{pMp}(P_2)'\cap pMp\subset L(\Gamma)'\cap pMp=\mathbb Cp$, the  claim follows from \cite[Lemma 2.4(3)]{DHI16}. Second, assume that $(1\otimes P_2)z\nprec_{M\bar\otimes M}M\bar\otimes (A\rtimes\Gamma_{\widehat{k_j}})$.
Since $Q_1^j z\subset ((1\otimes P_2)z)'\cap z (M\bar\otimes M) z$, \eqref{y1} implies that $((1\otimes P_2)z)'\cap z (M\bar\otimes M) z$ is strongly non-amenable relative to $M\bar\otimes (A\rtimes\Gamma_{\widehat {k_j}})$. Altogether, we can apply
Lemma \ref{lemma.trick.commultiplication.psi} to deduce that $(1\otimes \Psi^{k_j}) (z(M\bar\otimes L(\Lambda))z)$ is $({\rm id}\otimes{\rm id}\otimes \alpha^{k_j})$-rigid. Since $M\bar\otimes L(\Lambda)$ is a II$_1$ factor and $(1\otimes \Psi^{k_j})(1\otimes B)\subset 1\otimes B\otimes 1$, it follows that $(1\otimes \Psi^{k_j}) (M\bar\otimes M)$ is $({\rm id}\otimes{\rm id}\otimes \alpha^{k_j})$-rigid, which implies that $L(\Gamma_{k_j})$ is $\alpha^{k_j}$-rigid, contradiction. This completes the proof of the claim.
\hfill$\square$

Note that $Q_1$ and $(1\otimes P_2)q_1$ are commmuting subalgebras of $q_1 (M\bar\otimes M) q_1$ Thus, \eqref{y1} together with Lemma \ref{lemma.trick.commultiplication.psi} imply that for any $j\in S_1$ we have
\begin{equation}\label{y2}
    (1\otimes\Psi^{k_j})(\bigvee_{i\in S_1\setminus\{j\}}Q_1^iz \vee (1\otimes P_2)z) \text{ is } ({\rm id}\otimes{\rm id}\otimes \alpha^{k_j})\text{-rigid.}
\end{equation}

We now ready to prove the following.

{\bf Claim 3.} The map $S_1\ni j \to k_j\in \{1,\dots,n\}$ is injective.

{\it Proof of Claim 3.} Assume by contradiction that there exist two distinct elements $j_1,j_2\in S_1$ such that $k:=k_{j_1}=k_{j_2}$. Thus, $(S_1\setminus\{j_1\})\cup (S_1\setminus\{j_2\})=S_1$. Since the algebras $Q_1^jz, j\in S_1$, are commuting, we deduce from \eqref{y2} that
$(1\otimes\Psi^k)(Q_1z) \text{ is } ({\rm id}\otimes{\rm id}\otimes \alpha^{k})\text{-rigid.}$ As in the proof of Claim 1, we note that $zw_1\neq 0$.
Note also that  $Q_1zw_1=zw_1r_1\Delta(L(\Gamma_{S_1}))r_1$. By applying Proposition \ref{proposition.deformation.set.vaes} we obtain a non-zero projection $e_1\in (1\otimes\Psi^k)(\Delta(L(\Gamma)))'\cap M\bar\otimes M\bar\otimes L(\Gamma)$ such that 
\begin{equation}\label{y3}
(1\otimes\Psi^k)(\Delta(L(\Gamma_{S_1})))e_1 \text{ is } ({\rm id}\otimes{\rm id}\otimes \alpha^{k})\text{-rigid.}
\end{equation}
Since $z\in q_1(M\bar\otimes P_1)q_1$ and $M\bar\otimes P_1$ is a II$_1$ factor, one can check that \eqref{y2} implies 
\begin{equation}\label{ddextra}
\Psi^k(P_2) \text{ is } ({\rm id}\otimes \alpha^{k}) \text{ -rigid.}    
\end{equation}

Next, since $\Delta(L(\Gamma_{S_2}))\prec^s_{M\bar\otimes L(\Lambda)} M\bar\otimes P_2$, we obtain from \cite[Remark 2.2]{DHI16} that  
$$(1\otimes \Psi^k)(\Delta(L(\Gamma_{S_2})))\prec^s_{M\bar\otimes M\bar\otimes L(\Gamma_k)} M\bar\otimes \Psi^k(P_2).$$ 
Therefore, $(1\otimes \Psi^k)(\Delta(L(\Gamma_{S_2})))e_1\prec_{M\bar\otimes M\bar\otimes L(\Gamma_k)} M\bar\otimes \Psi(P_2)$, which implies by \eqref{ddextra} that there is a  projection $0\neq e_2\in (1\otimes\Psi^k)(\Delta(L(\Gamma_{S_2})))'\cap   (M\bar\otimes M\bar\otimes L(\Gamma_k))$ with $e_2\leq e_1$ such that 
\begin{equation}\label{y4}
(1\otimes\Psi^k)(\Delta(L(\Gamma_{S_2})))e_2 \text{ is } ({\rm id}\otimes{\rm id}\otimes \alpha^{k})\text{-rigid.}    
\end{equation}

Note that \eqref{y3} implies that $e_2(1\otimes\Psi^k)(\Delta(L(\Gamma_{S_1})))e_2 \text{ is } ({\rm id}\otimes{\rm id}\otimes \alpha^{k})\text{-rigid.}$ 
Together with \eqref{y4} and the fact that $\Psi^k(A)\subset A\otimes 1$, we deduce from Proposition \ref{proposition.deformation.set.vaes} that there exists a non-zero projection $e_3\in (1\otimes\Psi^k)(\Delta(M))'\cap M\bar\otimes M\bar\otimes L(\Gamma_k)$ such that 
\begin{equation}\label{y5}
(1\otimes\Psi^k)(\Delta(M))e_3 \text{ is } ({\rm id}\otimes{\rm id}\otimes \alpha^{k})\text{-rigid.}
\end{equation}
This implies that for any $\epsilon>0$, there exists $t_0>0$ such that for all $|t|\leq t_0$ and $g\in\Lambda$,
$$
\|(1\otimes \Psi^k)(v_g\otimes v_g)e_3 -({\rm id}\otimes{\rm id}\otimes \alpha_t^{k}) ((1\otimes \Psi^k)(v_g\otimes v_g)e_3)  \|_2\leq \epsilon,
$$
and, therefore,
$$
\|(1\otimes \Psi^k)(1\otimes v_g)e_3 -({\rm id}\otimes{\rm id}\otimes \alpha_t^{k}) ((1\otimes \Psi)(1\otimes v_g)e_3)  \|_2\leq \epsilon.
$$
Note that
$\Psi^k(B)\subset B\otimes 1$.
By applying Proposition \ref{proposition.deformation.set.vaes} we get that $\Psi^k(M)e_0$ is $({\rm id}\otimes \alpha^k)$-rigid for a projection $0\neq e_0\in \Psi^k(M)'\cap( M\bar\otimes L(\Gamma_k))$. Since $\Gamma_k$ is icc, we get $\Psi^k(M)'\cap( M\bar\otimes L(\Gamma_k))=\mathbb C1$. Thus, we obtain that
$L(\Gamma_k)$ is $\alpha^k$-rigid, contradiction.
\hfill$\square$

Denote $R_1=\{k_j| j\in S_1\}\subset \{1,\dots,n\}$. Claim 3 implies that $|S_1|=|R_1|$ while Claim 2 together with \cite[Lemma 2.8(2)]{DHI16} gives that $P_1\prec_M^s B\rtimes\Lambda_{\widehat{R_1}}$. In a similar way, there exists a subset $R_2\subset \{1,\dots,n\}$ with $|S_2|=|R_2|$ such that 
$P_2\prec_M^s A\rtimes\Gamma_{\widehat{R_2}}$. By using \cite[Proposition 4.4]{CD-AD21} we deduce that $L(\Gamma)\prec_M^s A\rtimes\Gamma_{\widehat{R_1}\cup \widehat{R_2}}$. Using \cite[Lemma 2.3]{BV12} we get that $M\prec_{M} A\rtimes\Gamma_{\widehat{R_1}\cup \widehat{R_2}}$, which implies that $\widehat{R_1}\cup \widehat{R_2}=\{1,\dots,n\}$.

Finally, we let $T_1=\widehat{R_1}$ and $T_2=\widehat{R_2}$. Since $S_1\sqcup S_2=\{1,\dots,n\}$ is a partition, it follows that $T_1\sqcup T_2=\{1,\dots,n\}$ is a partition as well. This ends the proof.
\hfill$\blacksquare$

\section{From unitary conjugacy of subalgebras to cohomologous cocycles}


In this section we prove Proposition \ref{proposition.unitary.cohomologous} which provides sufficient conditions at the von Neumann algebra level for untwisting the underlying cocycle of an orbit equivalence of irreducible actions. 

Throughout this section we will use the well known fact that if $\Gamma\car (X,\mu)$ and $\Lambda\car (Y,\nu)$ are free ergodic pmp actions such that there is a measure space isomorphism $\theta:X\to Y$ with $\theta(\Gamma x)=\Lambda\theta (x)$, for a.e. $x\in X$, then the induced isomorphism of von Neumann algebras $\pi: L^\infty(X) \to L^\infty(Y)$ given by $\pi(a)=a\circ\theta^{-1}$ extends to an isomorphism  $\pi: L^\infty(X)\rtimes\Gamma \to L^\infty(Y)\rtimes\Lambda$ satisfying $\pi(u_g)=v_{\theta\circ g \circ \theta^{-1}}$, for any $g\in\Gamma$. Here and throughout the section, we denote by $v_{\phi}\in \mathcal U(L^\infty(Y)\rtimes\Lambda)$ the associated unitary of $\phi\in [\mathcal R(\Lambda\car Y)]$; see \cite[Section 1.5.2]{AP10} for more details. 

\begin{proposition}\label{proposition.unitary.cohomologous}
Let $\Gamma=\Gamma_1\times\dots\times\Gamma_n\car (X,\mu)$ and $\Lambda=\Lambda_1\times\dots\times\Lambda_n\car (Y,\nu)$ be free, irreducible, pmp actions such that are orbit equivalent via a map $\theta:X\to Y$. Denote by $\pi:L^\infty(X)\rtimes\Gamma\to L^\infty(Y)\rtimes\Lambda$ the $*$-isomorphism associated to $\theta$ and let $c:\Gamma\times X\to\Lambda$ be the Zimmer cocycle associated to $\theta$.

If there exist $u_1,\dots,u_n\in\mathcal U(L^\infty(Y)\rtimes\Lambda)$ such that $\pi(L^\infty(X)\rtimes\Gamma_{\widehat i})= u_i(L^\infty(Y)\rtimes\Lambda_{\widehat i})u_i^*$, for any $i\in\{1,\dots,n\}$, then $c$ is cohomologous to a group isomorphism $\delta:\Gamma\to\Lambda$.

\end{proposition}

We first need the following elementary result. For completeness, we provide a proof.

\begin{lemma}\label{L1}
Let $\Gamma \overset{}{\car}  (X,\mu)$ and $\Lambda \overset{}{\car} (Y,\nu)$ be free, ergodic, pmp actions. For any $1\leq i\leq 2$, assume that there exist a $*$-isomophism $\pi_i: L^\infty(X)\rtimes\Gamma\to L^\infty(Y)\rtimes\Lambda$ such that $\pi_i (L^\infty(X))=L^\infty(Y)$, let $\theta_i:X\to Y$ be the measure space isomorphism defined by $\pi_i(a)=a\circ \theta_i^{-1}$, for any $a\in L^\infty(X)$, and let $c_i:\Gamma\times X\to\Lambda$ be the Zimmer cocycle associated to $\theta_i$.

If there exists $\omega\in \mathcal U(L^\infty(Y)\rtimes\Lambda)$ such that $\pi_2={\rm Ad}(\omega) \circ \pi_1$, then the cocycles $c_1$ and $c_2$ are cohomologous.
\end{lemma}

{\it Proof.} Since $\omega\in \mathcal N_{L^\infty(Y)\rtimes\Lambda}(L^\infty(Y))$, we can write $\omega=b v_{\varphi}$, for some $b\in \mathcal U(L^\infty(Y))$ and $\varphi\in [\mathcal R(\Lambda\car Y)]$ (see, for instance, \cite[Lemma 12.1.16]{AP10}). Take a measurable map $\psi: Y\to\Lambda$ such that $\varphi^{-1}(y)=\psi(y)y$, for almost every $y\in Y$.
For any $a\in L^\infty(X)$, we have $a\circ \theta_2^{-1}=\omega (a\circ \theta_1^{-1}) \omega^*=a\circ \theta_1^{-1}\circ \varphi^{-1}$. This shows that
$\theta_2=\varphi\circ \theta_1$. We will prove the lemma by showing that 
\begin{equation}\label{ff1}
    c_1(g,x)\psi(\theta_2(x))=\psi(\theta_2(gx))c_2(g,x), \text{ for all } g\in\Gamma \text{ and almost every }x\in X.
\end{equation}

To this end, fix an arbitrary  $g\in\Gamma$. Define $\tilde\psi=\psi\circ\theta_2$. Since for almost every $y\in Y$ and $i\in\{1,2\}$, we have  $(\theta_i\circ g \circ \theta_i^{-1})(y)=c_i(g,\theta_i^{-1}(y))y$, it follows that
\begin{equation*}\label{sss1}
\begin{array}{rcl}
c_1(g^{-1}, \theta_2^{-1}(y))\tilde \psi(\theta_2^{-1}(y))y &=&
c_1(g^{-1}, \theta_1^{-1}(\psi(y)y))\psi(y)y\\
&=&
(\theta_1\circ g \circ \theta_1^{-1})(\psi(y)y) )= (\theta_1\circ g\circ \theta_1^{-1}\circ \varphi^{-1})(y) \\
&=& (\varphi^{-1}\circ \theta_2\circ g \circ \theta_2^{-1})(y)=\psi ((\theta_2\circ g \circ \theta_2^{-1})(y)) (\theta_2\circ g \circ \theta_2^{-1})(y)\\
&=& \tilde\psi (g\theta_2^{-1}(y)) c_2(g,\theta_2^{-1}(y))y.
\end{array}    
\end{equation*}

Since $\Lambda\car Y$ is free, we obtain that \eqref{ff1}holds, thus proving the lemma.
\hfill$\blacksquare$

The following lemma is a particular case of \cite[Lemma 3.1]{HH21} and it goes back to \cite[Section 5]{MS02}. For the convenience of the reader, we provide a short proof for it using von Neumann algebras.


\begin{lemma}[\!\!\cite{HH21}]\label{L2}
Let $\Gamma=\Gamma_1\times\Gamma_2\overset{\sigma}{\car} (X,\mu)$ and $\Lambda=\Lambda_1\times\Lambda_2\overset{\rho}{\car} (Y,\nu)$ be free, ergodic, pmp actions with $\Gamma_1$ and $\Lambda_1$ acting ergodically. Assume that there exists a measure space isomorphism $\theta: X\to Y$ such that $\theta(\Gamma\cdot x)=\Lambda\cdot \theta(x)$ and $\theta(\Gamma_1\cdot x)=\Lambda_1\cdot \theta(x)$ for almost every $x\in X$. Let $c$ be the Zimmer cocycle associated to $\theta$.

Then there exists a group isomorphism $\delta_2:\Gamma_2\to \Lambda_2$ such that $c(g,x)\in \Lambda_1\delta_2(g_2)$ for every $g=(g_1,g_2)\in \Gamma$ and almost every $x\in X$.
\end{lemma}

{\it Proof.}  Denote by $\pi:L^\infty(X)\rtimes\Gamma\to L^\infty(Y)\rtimes\Lambda$ the $*$-isomorphism associated to $\theta$. For ease of notation, we suppress $\pi$.  Recall that for each $g\in\Gamma$ we can decompose
\begin{equation}\label{lh0}
u_g=\sum_{\lambda\in \Lambda}1_{Y_{g,\lambda}}v_\lambda, 
\end{equation}
where $Y_{g,\lambda}=\{y\in Y\; | c(g^{-1},\theta^{-1}(y))=\lambda^{-1} \},$ 
as $\lambda\in\Lambda$. By assumption, $c(g,x)\in\Lambda_1$, for any $g\in\Gamma_1$ and almost every $x\in X$.
Hence, we deduce $N:=L^\infty(X)\rtimes\Gamma_1=L^\infty(Y)\rtimes\Lambda_1$.

Next, we fix $g\in \Gamma_2$. Note that the actions $\sigma_{|\Gamma_2}$ and $\rho_{|\Lambda_2}$ extend in a natural way to actions on $N$. We can write $u_g=\sum_{\lambda\in\Lambda_2}b^g_\lambda v_\lambda$, with $b^g_\lambda\in N$, for all $\lambda\in\Lambda_2$. Note that for any $a\in N$ we have $b_\lambda^g\rho_\lambda(a)=\sigma_g(a)b_\lambda^g$, for all $\lambda\in \Lambda_2$. Thus, for any $\lambda\in\Lambda_2$, we get $(b_\lambda^g)^*b_\lambda^g\in N'\cap M=\mathbb C1$. Assume by contradiction that there exist $\lambda_1\neq \lambda_2\in\Lambda_2$ such that $b^g_{\lambda_1}$ and $b^g_{\lambda_2}$ are non-zero. Thus, there exist $\lambda_0\in\Lambda_2\setminus\{e\}$ and a unitary $c\in N$ such that $\rho_{\lambda_0}(a)c=ca$, for all $a\in N$. By writing $c=\sum_{\lambda\in\Lambda_1}c_\lambda v_\lambda$, we have $\rho_{\lambda_0\lambda^{-1}}(a)\rho_{\lambda^{-1}}(c_{\lambda})=a\rho_{\lambda^{-1}}(c_{\lambda})$, for all $a\in L^\infty(Y)$ and $\lambda\in\Lambda_1.$ Since $\lambda_0\lambda^{-1}$ acts freely, we get that $c=0$, contradiction. Thus, we have shown that there exist a map $\delta_2:\Gamma_2\to\Lambda_2$ and a unitary $b_g\in N$, as $g\in\Gamma_2$ satisfying
\begin{equation}\label{lh1}
u_g=b_gv_{\delta_2(g)}.
\end{equation}
One immediately obtains that $\delta_2:\Gamma_2\to\Lambda_2$ is a group homomorphism.
In a similar way, we can write $v_\lambda=\tilde b_{\lambda} u_{\eta_2(\lambda)}$ for some $\tilde b_\lambda\in N$ and a group homomorphism $\eta_2:\Lambda_2\to\Gamma_2$. It follows that $\eta_2\circ \delta_2={\rm Id}$ and $\delta_2\circ\eta_2={\rm Id}$, hence showing that $\delta_2$ is a group isomorphism.


Therefore, by combining \eqref{lh0} and \eqref{lh1}, we deduce that for any $g\in\Gamma_2$, we have that  $\mu(Y_{g,\lambda})=0$, for any $\lambda\notin \Lambda_1\delta_2(g)$. This  implies that $c(g,x)\in \Lambda_1\delta_2(g)$, for any $g\in\Gamma_2$. Finally, if $g=(g_1,g_2)\in\Gamma$, we get that $c(g,x)=c(g_1,g_2x)c(g_2,x)\in \Lambda_1\delta_2(g_2)$. This ends the proof of the lemma.
\hfill$\blacksquare$

\subsection{Proof of Proposition \ref{proposition.unitary.cohomologous}}
Fix an arbitrary $i\in\overline{1,n}$. By Lemma \ref{L1} we get that the underlying Zimmer cocycle $c_i:\Gamma\times X\to\Lambda$ of the orbit equivalence $\theta_i:X\to Y$ associated to $\pi_i:= {\rm Ad}(u_i^*)\circ\pi:L^\infty(X)\rtimes\Gamma\to L^\infty(Y)\rtimes\Lambda$ is cohomologous to $c$. 
Hence, there is a map $\varphi_i:X\to\Lambda$ such that $c(g,x)=\varphi_i(gx)^{-1}c_i(g,x)\varphi_i(x)$, for all $g\in\Gamma$ and almost every $x\in X$.
Note that $\Gamma_{\widehat i}\car X$ and $\Lambda_{\widehat i}\car Y$ are ergodic. Since $\pi_i(L^\infty(X)\rtimes\Gamma_{\widehat i})=L^\infty(Y)\rtimes\Lambda_{\widehat i}$, we get that $\theta_i(\Gamma_{\widehat i}\cdot x)=\Lambda_{\widehat i}\cdot \theta(x)$, for almost every $x\in X$, and hence, we obtain from Lemma \ref{L2} that there is a group isomorphism $\delta_i:\Gamma_i\to \Lambda_i$ such that $c_i(g,x)\in 
\Lambda_{\widehat i}\delta_i(g_i)$, for every $g=(g_{\widehat i},g_i)\in \Gamma=\Gamma_{\widehat i}\times\Gamma_i$ and almost every $x\in X$. 

Next, since $\Lambda=\Lambda_1\times\dots\times\Lambda_n$ we decompose $\varphi_i=\varphi_i^1\dots \varphi_i^n$ and $c_i=c_i^1 \dots c_i^n$ where $\varphi_i^j$ and $c_i^j$ are valued to $\Lambda_j$, for any $j\in\overline{1,n}$. By letting $\varphi=\varphi_1^1\dots\varphi_n^n:X\to\Lambda$ and $\tilde c:\Gamma\times X\to\Lambda$ defined by $\tilde c(g,x)=\varphi(gx)c(g,x)\varphi(x)^{-1}$, we get
$\tilde c(g,x)=\phi(gx) \phi_i(gx)^{-1} c_i(g,x) \phi_i(x)  \phi(x)^{-1}$, for all $i\in\overline{1,n}$, $g\in\Gamma$ and almost every $x\in X$. Consequently, we obtain that for every $g=(g_{\widehat i},g_i)\in \Gamma=\Gamma_{\widehat i}\times\Gamma_i$ and almost every $x\in X$, we have
$\tilde c^i(g,x)=c_i^i(g,x)=\delta_i(g_i);$    
here, we denoted by $\tilde c=\tilde c^1 \dots \tilde c^n$ the decomposition along $\Lambda=\Lambda_1\times\dots\times\Lambda_n$.
We define the group isomorphism  $\delta:\Gamma\to\Lambda$ by letting $\delta(g_1\dots g_n)=\delta_1(g_1)\dots \delta_n(g_n)$, for all $g_1\in\Gamma_1,\dots,g_n\in\Gamma_n$. This shows that $\tilde c (g,x)=\delta(g)$, for all $g\in\Gamma$ and almost every $x\in X$, which entails to $c$ is cohomologous to the group isomorphism $\delta.$
\hfill$\blacksquare$

\section{Proofs of Theorems \ref{B_0}, \ref{A} and Corollary \ref{Main.new.prime}}

In this section we prove the first three main results stated in the introduction. 
Towards this we first prove an abstract version of \cite[Lemma 5.10]{DHI16} in the sense that we only require the lack of property Gamma instead of relative solidity assumptions. In order to properly state and prove the result, we assume the terminology introduced in Notation \ref{notation}.

\begin{lemma}\label{Step3}
Let $L(\Lambda)=P_1\bar\otimes P_2$ be a tensor product decomposition into II$_1$ factors.

Assume there exist two partitions $T_1\sqcup T_2=S_1\sqcup S_2=\{1,\dots,n\}$ such that for any $i\in \{1,2\}$ we have  $P_i\prec^s_{M} A\rtimes \Gamma_{T_i}$ and $\Delta(A\rtimes\Gamma_{S_i})\prec^s_{M\bar\otimes L(\Lambda)} M\bar\otimes P_i$.

If $L(\Lambda)$ does not have property Gamma, then $T_i=S_i$, for any $i\in\{1,2\}$. Moreover, there exist 
subgroups $\Sigma_1,\Sigma_2<\Lambda$ such that for all $i\in \{1,2\}$ we have:
\begin{enumerate}
    \item $B\rtimes\Sigma_i\prec^s_M A\rtimes \Gamma_{T_i}$ and $A\rtimes\Gamma_{S_i}\prec^s_M B\rtimes\Sigma_i$,
    \item $P_i\prec^s_{L(\Lambda)} L(\Sigma_i)$ and $L(\Sigma_i)\prec_{L(\Lambda)}^s P_i$.
\end{enumerate}
\end{lemma}

{\it Proof.} (1) The assumption $\Delta(A\rtimes\Gamma_{S_1})\prec^s_{M\bar\otimes L(\Lambda)} M\bar\otimes P_1$ implies from \cite[Theorem 4.1]{DHI16} (see also \cite[Theorem 3.1]{Io11} and \cite[Theorem 3.3]{CdSS15}) that there exists a decreasing sequence of subgroups
$\Omega_k<\Lambda$, $k\ge 1$, such that $A\rtimes\Gamma_{S_1}\prec_M^s B\rtimes\Omega_k$, for any $k\ge 1$ and $P_2\prec_{L(\Lambda)} L(\cup_{k\ge 1}C_{\Lambda}(\Omega_k))$. Using Proposition \ref{proposition.full.sequence}, there is $k\ge 1$ such that $P_2\prec_{L(\Lambda)} L(\Omega_k)'\cap L(\Lambda)$ and using
\cite[Lemma 3.5]{Va08} we further derive that $L(\Omega_k)\prec_{L(\Lambda)} P_1$. By letting $\Sigma_1=\Omega_k$, we get
\begin{equation}\label{a1}
    A\rtimes\Gamma_{S_1}\prec_M^s B\rtimes\Sigma_1 \text{ and } L(\Sigma_1)\prec_{L(\Lambda)} P_1.
\end{equation}

We continue by showing that $B\rtimes\Sigma_1\prec^s_M A\rtimes \Gamma_{T_1}$.
By applying \cite[Lemma 2.4]{DHI16}, we get from \eqref{a1} a non-zero projection $e\in L(\Sigma_1)'\cap L(\Lambda)$ such that $L(\Sigma_1)e\prec^s_{L(\Lambda)} P_1$. 
 Since $P_1\prec^s_{M} A\rtimes \Gamma_{T_1}$, we obtain from \cite[Lemma 3.7]{Va08} that $L(\Sigma_1)e\prec^s_M A\rtimes \Gamma_{T_1}$. By applying Lemma \ref{lemma.joint.control} there exists a projection $f\in (B\rtimes\Sigma_1)'\cap pMp\subset B$ with $f\ge e$ such that $(B\rtimes\Sigma_1)f\prec^s_M A\rtimes \Gamma_{T_1}$. Since $f\in B$, $e\in L(\Lambda)$ and $f\ge e$, we deduce that $f=1$. Thus, $B\rtimes\Sigma_1\prec^s_M A\rtimes \Gamma_{T_1}$. Similarly, there exists a subgroup $\Sigma_2<\Lambda$ satisfying conclusion (1).
 
 (2) This follows verbatim the proofs of Claims 2 and 3 from \cite[Lemma 5.10]{DHI16}.
\hfill$\blacksquare$




\subsection{Proof of Theorem \ref{B_0}}
Assume Notation \ref{notation}. Fix an arbitrary $i\in\overline{1,m}$. By applying Theorem \ref{Step1} there is a partition $S^i_1\sqcup S^i_2=\overline{1,n}$ such that $\Delta(A\rtimes\Gamma_{S_1^i})\prec_{M\bar\otimes L(\Lambda)} M\bar\otimes L(\Lambda_{\widehat i})$ and $\Delta(A\rtimes\Gamma_{S_2^i})\prec_{M\bar\otimes L(\Lambda)} M\bar\otimes L(\Lambda_{i})$. Standard arguments imply that
\begin{equation}\label{tts1}
A\rtimes\Gamma_{S_1^i}\prec_{M}^s B\rtimes\Lambda_{\widehat i} \text{ and }    A\rtimes\Gamma_{S_2^i}\prec_{M}^s B\rtimes\Lambda_{ i}.
\end{equation}
Hence, Theorem \ref{Step2} combined with \cite[Lemma 2.3]{BV12} gives that 
there is a partition $S^i_1\sqcup S^i_2=\overline{1,n}$ such that 
\begin{equation}\label{tts2}
    B\rtimes\Lambda_{\widehat i}\prec^s_{M} A\rtimes\Gamma_{T_1^i} \text{ and } B\rtimes\Lambda_{ i}\prec^s_{M} A\rtimes\Gamma_{T_2^i}.
\end{equation}
By applying \cite[Lemma 3.7]{Va08} we derive that $S_1^i=T_1^i$, $S_2^i=T_2^i$.
Consequently, by using relations \eqref{tts1} and \eqref{tts2},
\cite[Proposition 3.1]{DHI16} implies that $\Gamma_{T_1^i}$ and $\Lambda_{\widehat i}$ are measure equivalent and  $\Gamma_{T_2 ^i}$ and $\Lambda_{i}$ are measure equivalent as well, for any $i\in\overline{1,m}$. The conclusion now follows by a simple induction argument.
\hfill$\blacksquare$

\subsection{Proof of Theorem \ref{Theorem.ME}}
We first obtain the following classification of tensor product decompositions in the spirit of \cite[Theorem C]{DHI16}.   Theorem \ref{Theorem.ME} will then follow by applying this result together with an induction argument.


\begin{theorem}\label{theorem.classification.tensor.product}
Let $\Gamma$ and  $\Lambda$ be countable icc groups that are measure equivalent. Assume $L(\Lambda)=P_1\bar\otimes P_2$ and $\Gamma=\Gamma_1\times\dots\times\Gamma_n$ is a product into icc groups such that $L(\Gamma_i)$
belongs to $\bm{\mathscr M}$ for any $i\in\{1,\dots,n\}$.

Then there exist a direct product decomposition $\Lambda=\Lambda_1\times\Lambda_2$, a partition $T_1\sqcup T_2=\{1,\dots,n\}$, a decomposition $L(\Lambda)=P_1^{t_1}\bar\otimes P_2^{t_2}$, for some $t_1,t_2>0$ with $t_1t_2=1$, and a unitary $u\in L(\Lambda)$ such that:
\begin{enumerate}
    \item $P_1^{t_1}=uL(\Lambda_1)u^*$ and $P_2^{t_2}=uL(\Lambda_2)u^*$. 
    \item $\Lambda_1$ is measure equivalent to $\times_{j\in T_1}\Gamma_j$ and $\Lambda_2$ is measure equivalent to $\times_{j\in T_2}\Gamma_j$.
\end{enumerate}
\end{theorem}

{\it Proof.} For the proof, we assume Notation \ref{notation}.
Using Proposition \ref{property.gamma}, we get that $L(\Gamma)$ does not have property Gamma. 
Next, by applying Theorem \ref{Step1}, Theorem \ref{Step2} and Lemma \ref{Step3}, we obtain a partition $T_1\sqcup T_2=\{1,\dots,n\}$ and some subgroups $\Sigma_1,\Sigma_2<\Lambda$ such that for all $i\in \{1,2\}$ we have:
\begin{enumerate}
    \item $B\rtimes\Sigma_i\prec^s_M A\rtimes \Gamma_{T_i}$ and $A\rtimes\Gamma_{S_i}\prec^s_M B\rtimes\Sigma_i$,
    \item $P_i\prec^s_{L(\Lambda)} L(\Sigma_i)$ and $L(\Sigma_i)\prec_{L(\Lambda)}^s P_i$.
\end{enumerate}

Part (2) together with \cite[Theorem 6.1]{DHI16} give a product decompomsition $\Lambda=\Lambda_1\times\Lambda_2$,
a decomposition $L(\Lambda)=P_1^{t_1}\bar\otimes P_2^{t_2}$, for some $t_1,t_2>0$ with $t_1t_2=1$, and a unitary $u\in L(\Lambda)$ such that
$P_1^{t_1}=uL(\Lambda_1)u^*$ and $P_2^{t_2}=uL(\Lambda_2)u^*$. In addition, we have that $\Lambda_i$ is measure equivalent to $\Sigma_i$, for any $i\in\{1,2\}$.

Part (1) together with \cite[Proposition 3.1]{DHI16} implies that for any $i\in\{1,2\}$, $\Gamma_{T_i}$ is measure equivalent to $\Sigma_i,$ and hence, to $\Lambda_i$.
\hfill$\blacksquare$




\subsection{Proof of Corollary \ref{Main.new.prime}}
Assume first that there exists a partition $J_1\sqcup J_2=\{1,\dots,n\}$ for which $K={\rm Fix}_K(J_1)\times {\rm Fix}_K(J_2)$. By letting $G_1= (\times_{i\in J_1} \Gamma)\rtimes  {\rm Fix}_K(J_2)$ and $G_2=(\times_{i\in J_2}\Gamma)\rtimes  {\rm Fix}_K(J_1)$, it follows that $G= G_1 \times  G_2$. This clearly shows that $L(G)$ is not prime.

For proving the other implication, assume that $L(G)=P_1\bar\otimes P_2$ can be written as a tensor product of diffuse factors. Using Theorem \ref{A} and its proof, it follows that there exist a direct product decomposition $G=G_1\times G_2$ into infinite groups 
and a partition $J_1\sqcup J_2=\overline{1,n}$ such that $L(G_i)\prec_{L(G)} L(\Gamma_{J_i})$, for any $i\in\overline{1,2}.$
By \cite[Lemma 2.2]{CI17} we get a finite index subgroup $G_i^0<G_i$ such that $G_i^0< \Gamma_{J_i}$, for any $i\in\overline{1,2}$. 
By passing to relative commutants, we get that $\Gamma_{J_2}< G_2$ since $G_1$ is icc. By passing again to relative commutants, we deduce that $G_1< C_{G}(\Gamma_{J_2})=(\times_{i\in J_1} \Gamma)\rtimes  {\rm Fix}_K(J_2)$.  Similarly, we get $G_2<(\times_{i\in J_2} \Gamma)\rtimes  {\rm Fix}_K(J_1)$. This proves $K={\rm Fix}_K(J_1)\times {\rm Fix}_K(J_2)$ which ends the proof.
\hfill$\blacksquare$

\section{Proof of Theorem \ref{C}}

\subsection{OE rigidity for irreducible actions}

An important ingredient for proving Theorem \ref{C} is the following OE rigidity result for irreducible actions of product group that belong to $\bm{\mathscr{M}}.$

\begin{theorem}\label{B}
Let $\Gamma=\Gamma_1\times\dots\times\Gamma_n$ be a product of $n\ge 2$ groups that belong to $\bm{\mathscr{M}}$. Let $\Lambda=\Lambda_1\times\dots\times \Lambda_m$ be a product of $m\ge 2$ infinite icc groups. Assume $\Gamma\car (X,\mu)$ and $\Lambda\car (Y,\nu)$ are orbit equivalent  free, irreducible, pmp actions.

If $m\ge n$, then $m=n$ and $\Gamma\car (X,\mu)$ and $\Lambda\car (Y,\nu)$ are conjugate.
\end{theorem}

By assumption we have the identification $M=L^\infty(X)\rtimes\Gamma=L^\infty(Y)\rtimes\Lambda$ with $A:=L^\infty(X)=L^\infty(Y)$. Let $\Delta:M\to M\bar\otimes L(\Lambda)$ be the $*$-homomorphism given by $\Delta(av_\lambda)=av_\lambda\otimes v_\lambda$, for all $a\in A$ and $\lambda\in\Lambda$.
Fix $i\in\overline{1,m}$. By applying Theorem \ref{Step1} there is a partition $S^i_1\sqcup S^i_2=\overline{1,n}$ such that $\Delta(A\rtimes\Gamma_{S_1^i})\prec_{M\bar\otimes L(\Lambda)} M\bar\otimes L(\Lambda_{\widehat i})$ and $\Delta(A\rtimes\Gamma_{S_2^i})\prec_{M\bar\otimes L(\Lambda)} M\bar\otimes L(\Lambda_{i})$. Standard arguments imply that
\begin{equation}\label{ts1}
A\rtimes\Gamma_{S_1^i}\prec_{M}^s A\rtimes\Lambda_{\widehat i} \text{ and }    A\rtimes\Gamma_{S_2^i}\prec_{M}^s A\rtimes\Lambda_{ i}.
\end{equation}
Theorem \ref{Step2} combined with \cite[Lemma 2.3]{BV12} gives that 
there is $k_i\in \overline{1,n}$ such that 
\begin{equation}\label{ts2}
    A\rtimes\Lambda_{\widehat i}\prec^s_{M}A\rtimes\Gamma_{\widehat{k_i}} \text{ and } A\rtimes\Lambda_{ i}\prec^s_{M} A\rtimes\Gamma_{{k_i}}.
\end{equation}
By applying \cite[Lemma 3.7]{Va08} we derive that $m=n$, $S_1^i= \widehat{k_i}$ and $S_2^i=\{k_i\}$.
By using \cite[Lemma 8.4]{IPP05} there is $u_i\in \mathcal U(M)$ such that $u_i (A\rtimes\Lambda_{\widehat i}) u_i^*=A\rtimes \Gamma_{\widehat {k_i}}$. Thus, we can apply Proposition \ref{proposition.unitary.cohomologous} and derive that the Zimmer cocycle associated to the orbit equivalence between $\Gamma\car (X,\mu)$ and $\Lambda\car (Y,\nu)$ is cohomologous to a group isomorphism. Hence, by applying \cite[Lemma 4.7]{Va06} we get that $\Gamma\car (X,\mu)$ and $\Lambda\car (Y,\nu)$ are conjugate.
\hfill$\blacksquare$

\subsection{Strongly cocycle rigidity}
We start this subsection by recording the following particular case of \cite[Theorem 7.1]{HH21} which is inspired by several works \cite{Fu99,MS02,Ki08}. For properly formulating the result we introduce the following definition (see also \cite[Section 7]{HH21}). We say that a product group $\Gamma=\Gamma_1\times\dots\times\Gamma_n$ is {\it strongly cocycle rigid} if given any two free, irreducible, pmp actions $\Gamma\car (X,\mu)$ and $\Gamma\car (Y,\nu)$ that are orbit equivalent, the underlying Zimmer cocycle is cohomologous to a group isomorphism.

\begin{theorem}[\!\!\cite{HH21}]\label{HH21}
Let $\Gamma=\Gamma_1\times\dots\times\Gamma_n$ be an icc strongly cocycle rigid group. Assume $\Gamma\car(X,\mu)$ is a free, irreducible, pmp action that is orbit equivalent to a free, mildly mixing, pmp action $\Lambda\car (Y,\nu)$.

Then $\Gamma\car(X,\mu)$ and $\Lambda\car (Y,\nu)$ are virtually conjugate.
\end{theorem}

\begin{corollary}\label{corollary.strongly.cocycle.rigid}
If $\Gamma_1,\dots,\Gamma_n$ are countable groups with $L(\Gamma_i)\in  \bm{\mathscr M}$, for any $i\in\{1,\dots,n\}$, then $\Gamma_1\times\dots\times\Gamma_n$ is strongly cocycle rigid. 
\end{corollary}

{\it Proof.} Denote $\Gamma=\Gamma_1\times\dots\times\Gamma_n$ and let $\Gamma\car (X,\mu)$ and $\Gamma\car (Y,\nu)$ be free, irreducible, pmp actions that are orbit equivalent. 
The proof of Theorem \ref{B} gives, in particular, that
the underlying Zimmer cocycle is cohomologous to a group isomorphism. 
\hfill$\blacksquare$

\subsection{Proof of Theorem \ref{C}}
This is a direct consequence of Corollary \ref{corollary.strongly.cocycle.rigid} and Theorem \ref{HH21}.
\hfill$\blacksquare$

\section{Proofs of Theorem \ref{main.theorem.inneramenable} and Corollary \ref{mcor.inner.amenable}}

\subsection{Proof of Theorem \ref{D}}
Let $\{u_g\}_{g\in\Gamma}$ be the canonical unitaries that generate $L(\Gamma)$.
Denote $\mathcal M=M\bar\otimes L(\Gamma)$, $\tilde{\mathcal M}=\tilde M\bar\otimes L(\Gamma)$ and $\hat\alpha_t=\alpha_t\otimes {\rm id}\in {\rm Aut}(\tilde{\mathcal M})$.
Note that the $*$-homomorphism $\Delta:L(\Gamma)\to L(\Gamma)\bar\otimes L(\Gamma)$   defined by $\Delta(u_g)=u_g\otimes u_g$, as $g\in\Gamma$ \cite{PV09},  naturally extends to a map $\Delta:\ell^2(\Gamma)\to \ell^2(\Gamma)\otimes \ell^2(\Gamma)$. By denoting $\hat \xi=\Delta(\xi)$, for any $\xi\in \ell^2(\Gamma)$, it follows that if $\xi=\sum_{g\in\Gamma}\xi_g u_g\in\ell^2(\Gamma)$ and $t\in \mathbb R$, then
\begin{equation}\label{t1}
    \| \hat\alpha_t (\hat \xi)-\hat\xi \|_2^2=\sum_{g\in\Gamma}|\xi_g|^2 \| \alpha_t(u_g)-u_g \|_2^2.
\end{equation}
Since $\Gamma$ is inner amenable, there exists a sequence $(\xi_n)_{n\ge 1}\subset \ell^2(\Gamma)$ of unit vectors satisfying $\|u_g \xi_n -\xi_n u_g\|_2\to 0,$ for all $g\in\Gamma$ and $\xi_n(g)\to 0$, for any $g\in\Gamma$.
Let $\omega$ be a free ultrafilter on $\mathbb N$. The remaining part of the proof is divided between two claims.

{\bf Claim 1.} $\lim_{t\to 0}  (\lim_{n\to\omega} \| \hat\alpha_t(\hat \xi_n) -\hat\xi_n \| )=0$.

{\it Proof of Claim 1.} We define the unitary representations $\pi:\Gamma\to \mathcal U(L^2(\tilde M)\ominus L^2(M))$ by $\pi_g(\xi)=u_g\xi u_g^*$, for all $g\in\Gamma, \xi\in L^2(\tilde M)\ominus L^2(M)$ and $d:\Gamma\to \mathcal U(\ell^2(\Gamma))$ by $d_g(x)=u_g xu_g^*$, for all $g\in\Gamma, x\in \ell^2(\Gamma)$.  Since $L^2(\tilde M)\ominus L^2(M)$ is weakly contained in the coarse bimodule $L^2(M)\otimes L^2(M)$ as $M$-bimodules, we derive that $L^2(\tilde M)\ominus L^2(M)$ is weakly contained in the coarse bimodule $L^2(L(\Gamma))\otimes L^2(L(\Gamma))$ as $L(\Gamma)$-bimodules. Therefore,
 $\pi$ is weakly contained in the left regular representation $\lambda_{\Gamma}$. 
Consequently, by applying \cite[Corollary E.2.6]{BHV07} we derive that $\pi\otimes d$ is weakly contained in $\lambda_{\Gamma}$. Note that $\hat\pi:=\pi\otimes d :\Gamma\to \mathcal U(L^2(\tilde {\mathcal M})\ominus L^2(\mathcal M))$ is defined by $\hat\pi_g(\eta)=\hat{u}_g\eta \hat{u}_g^*$ for all $g\in\Gamma,\eta\in L^2(\tilde {\mathcal M})\ominus L^2(\mathcal M)$. 
Since $\Gamma$ is non-amenable, it follows that the trivial representation $1_{\Gamma}$ is not weakly contained in $\hat\pi$. This implies that for any $\epsilon>0$, there exist $\delta>0$ and a finite set $F\subset\Gamma$ satisfying that for any unit vector $\eta\in L^2(\tilde{\mathcal M})$ for which $\|\hat\pi_g(\eta)-\eta\|_2\leq\delta$, as $g\in F$, we have 
\begin{equation}\label{t0}
\| \eta - E_{\mathcal M}(\eta) \|_2\leq\epsilon.    
\end{equation}
Since $\tau (\hat\alpha_t(\hat u_g)\hat u_h^*)=\tau (\alpha_t( u_g) u_h^*) \delta_{g,h}$, we obtain
 that $E_{\Delta(L(\Gamma))}(\hat\alpha_t(\hat u_g))=\tau (\alpha_t(u_g)u_g^*)\hat u_g$, for any $g\in\Gamma$. 
  This  implies that 
 for all $g\in\Gamma$ and $\xi\in\ell^2(\Gamma)$, we have 
\begin{equation}\label{t2}
\| \hat\alpha_t(\hat u_g) \hat\xi -\hat u_g\hat \xi \|_2=\| \alpha_t( u_g)  - u_g \|_2 \|\xi\|_2    \text{ and }
\|\hat\xi \hat\alpha_t(\hat u_g)  -   \hat\xi\hat u_g \|_2=\| \alpha_t(u_g)  -u_g \|_2 \|\xi\|_2.  
\end{equation}
Let $t_0>0$ such that $\| \alpha_t(u_g)-u_g  \|_2\leq \delta/4$, for all $|t|<t_0$ and $g\in F$. Take also $n_0\in \mathbb N$ such that $\| u_g\xi_n -\xi_n u_g \|_2\leq \delta/2$, for all $g\in F$ and $n\ge n_0$. Together with \eqref{t2} we obtain
\begin{equation}\label{t3}
\begin{array}{rcl}
\| \hat\alpha_t(\hat u_g) \hat\xi_n -  \hat \xi_n \hat\alpha_t(\hat u_g) \|_2 &\leq&  \| \hat\alpha_t(\hat u_g) \hat\xi_n -\hat u_g\hat \xi_n \|_2 +
\| u_g\xi_n -\xi_n u_g \|_2
+\|\hat\xi_n\hat u_g-  \hat\xi_n \hat\alpha_t(\hat u_g)     \|_2\\
&\leq&  \delta/2 +\delta/4+\delta/2=\delta,
\end{array}
\end{equation}
for all $g\in F, n\ge n_0$ and $|t|\leq t_0$. By applying $\hat\alpha_{-t}$ in \eqref{t3} and by replacing $t$ by $-t$, we get that  $\| \hat\alpha_t(\hat\xi_n) \hat u_g -  \hat u_g \hat\alpha_t(\hat \xi_n) \|_2\leq \delta$, for all $g\in F, n\ge n_0$ and $|t|\leq t_0$. Using \eqref{t0}, we get $\| \hat\alpha_t(\hat\xi_n) - E_{\mathcal M}(\hat\alpha_t(\hat\xi_n))  \|_2\leq \epsilon$, and by using Popa's transversality property, see \cite[Lemma 2.1]{Po06a}, we further derive that $\| \hat\alpha_{2t}(\hat\xi_n))-\hat\xi_n   \|_2\leq 2\epsilon$,
for all $n\ge n_0$ and $|t|\leq t_0$. This ends the proof of the claim.
\hfill$\square$

For all $t\in\mathbb R$ and $r>0$, we denote $B_r^t=\{ g\in\Gamma | \; \|\alpha_t(u_g)-u_g \|_2\leq r \}$. We are now ready to prove the following claim.

{\bf Claim 2.} $\lim_{t\to 0}(\sup_{g\in\Gamma} \|\alpha_t(u_g)-u_g  \|_2)=0$.

{\it Proof of Claim 2.}
To this end, fix some arbitrary $\epsilon>0$. Let $t_1>0$ and $n_1\in \mathbb N$ such that $\| \hat\alpha_t(\hat \xi_n) -\hat\xi_n  \|_2\leq \epsilon/4 $, for all $|t|\leq t_1$ and $n\ge n_1$. Fix $g\in\Gamma$ and $|t|\leq t_1$.
We continue by showing that there exists an unbounded sequence $(k_n)_n \subset\Gamma$ such that $k_n, g k_n g^{-1}\in B_{\frac{\epsilon}{2}}^t$, for any $n\ge 1$. Since $\|\hat u_g \hat\xi_n -\hat\xi_n \hat u_g\|_2\to 0$, we get that there exists $n_2\in\mathbb N$ such that 
$$
\| \hat\alpha_t(\hat \xi_n) -\hat\xi_n  \|_2^2 + \| \hat\alpha_t(\hat u_g \hat \xi_n \hat u_g^*) - \hat u_g \hat\xi_n \hat u_g^*  \|_2^2 \leq \epsilon^2/4, \text{ for any }  n\ge n_2.
$$
By writing $\xi_n=\sum_{g\in\Gamma} \xi_{n,g} u_g\in \ell^2(\Gamma)$ and using
\eqref{t1} we obtain that
$$
\sum_{h\in\Gamma} |\xi_{n,h}|^2 (  \|\alpha_t(u_h)-u_h \|_2^2 + \|\alpha_t(u_{ghg^{-1}})-u_{ghg^{-1}} \|_2^2  )\leq \epsilon^2/4, \text{ for any }  n\ge n_2.
$$
For any $n\ge n_2$, since $\sum_{h\in\Gamma} |\xi_{n,h}|^2=1$, there exists $k_n\in\Gamma$ with $\xi_{n,k_n}\neq 0$ such that $  \|\alpha_t(u_{k_n})-u_{k_n} \|_2\leq\epsilon/2$ and $\|\alpha_t(u_{gk_ng^{-1}})-u_{gk_ng^{-1}} \|_2  \leq \epsilon/2$. Since $\xi_{n,h}\to 0$, for any $h\in\Gamma$, it follows that $k_n$ can be chosen such that $k_n\to\infty$. 

Next, we note that that for any $n\ge n_2$, we have
\begin{equation}\label{t4}
    \| \alpha_t (u_{gk_ng^{-1}})u_g - u_g \alpha_t(u_{k_n})  \|_2\leq \|\alpha_t(u_{k_n})-u_{k_n} \|_2 + \|\alpha_t(u_{gk_ng^{-1}})  - u_{gk_n g^{-1}} \|_2\leq \epsilon.
\end{equation}
By letting $e:L^2 (\tilde M)\to L^2(M)$ be the orthogonal projection, we have $v_{g,t}:=\alpha_t(u_g)-e(\alpha_t(u_g))\in L^2(\tilde M)\ominus L^2(M)$.
By applying $\alpha_{-t}$ to \eqref{t4} and by projecting onto $L^2(\tilde M)\ominus L^2(M)$, we get
\begin{equation}\label{t5}
\| u_{gk_n g^{-1}} v_{-t,g} -v_{-t,g} u_{k_n}   \|_2\leq \epsilon, \text{ for all } n\ge n_2.   
\end{equation}
Since the $M$-bimodule $L(\tilde M)\ominus L^2(M)$ is mixing
 and $k_n\to\infty$, we obtain that 
 \begin{equation}\label{t6}
 \lim_{n\to\infty}\langle v_{-t,g}u_{k_n},u_{gk_ng^{-1}} v_{-t,g}   \rangle=0.    
 \end{equation}
 
 By combining \eqref{t5} and \eqref{t6}, it follows that $\|v_{-t,g}\|_2\leq \epsilon/\sqrt{2}$. By using once again Popa's transversality property, see \cite[Lemma 2.1]{Po06a}, we obtain $\| \alpha_{-t}(u_g)-u_g \|\leq \epsilon \sqrt{2}$. Since $t$ was arbitrary chosen such that $|t|\leq t_1$ and $g\in\Gamma$ arbitrary, we get that   $\lim_{t\to 0}(\sup_{g\in\Gamma} \|\alpha_t(u_g)-u_g  \|_2)=0$.
\hfill$\square$

Standard arguments imply now the conclusion. 
\hfill$\blacksquare$

\subsection{Proof of Corollary \ref{mcor.inner.amenable}}
The proof follows directly from Theorem \ref{main.theorem.inneramenable}.
\hfill$\blacksquare$

\subsection{Consequence to Kurosh-type rigidty results}
We conclude our paper with the following rigidity result for tracial free product factors arising from non-amenable inner amenable groups.

\begin{corollary}\label{corollary.kurosh}
Let $M=L(\Gamma_1)*\dots*L(\Gamma_m)=L(\Lambda_1)*\dots*L(\Lambda_n)$, where all the groups $\Gamma_i$ and $\Lambda_j$ are non-amenable inner amenable icc groups.

Then $m=n$, and after a permutation of indices, $L(\Gamma_i)$ is unitarily conjugate to $L(\Lambda_i)$, for any $i\in\overline{1,n}.$
\end{corollary}

{\it Proof.}
Fix an arbitrary $i\in\overline{1,m}$. 
By decomposing $M=L(\Gamma_1*\dots*\Gamma_{n-1})*L(\Gamma_n)$, we note that $M$ belongs to $\bm{\mathscr M}$ and let $(\tilde M, (\alpha_t)_{t\in\mathbb R})$ be the associated  s-malleable deformation of $M$. Since $\Gamma_i$ is non-amenable inner amenable group,  Theorem \ref{D} implies that $L(\Gamma_i)$ is $\alpha$-rigid. By applying the main technical result of \cite{IPP05} (see also \cite[Theorem 2.11]{Io12a}), we get that $L(\Gamma_i)\prec_M L(\Gamma_1*\dots*\Gamma_{n-1})$ or $L(\Gamma_i)\prec_M L(\Gamma_n)$. 
By assuming the latter holds, there exist projections $p\in L(\Gamma_i), q\in L(\Gamma_n)$, a non-zero partial isometry $v\in qMp$ and a $*$-homomorphism $\theta: pL(\Gamma_i)p\to qL(\Gamma_n)q$ satisfying $\theta(x)v=vx$, for all $x\in pL(\Gamma_i)p$. Note that \cite[Theorem 1.2.1]{IPP05} gives that $vv^*\in L(\Gamma_n)$, and hence, $v L(\Gamma_i)v^*\subset L(\Gamma_n)$.
Note that since $L(\Gamma_i)'\cap M=\mathbb C 1$ and $v^*v\in p(L(\Gamma_i)'\cap M) p$, we get that $v^*v=p$.
By letting $u$ be a unitary that extends $v$, we derive that $u p L(\Gamma_i) pu^*\subset L(\Gamma_n)$. Since $L(\Gamma_n)$ is a factor, after passing to a new unitary $u$, one can replace $p$ by its central support in $L(\Gamma_i)$; therefore, we obtain that $u L(\Gamma_i)u^*\subset L(\Gamma_n)$. Similarly, if $L(\Gamma_i)\prec_M L(\Gamma_1*\dots*\Gamma_{n-1})$ holds, we obtain a unitary $u\in M$ such that $u L(\Gamma_i)u^*\subset L(\Gamma_1*\dots*\Gamma_{n-1})$. By repeating this argument finitely many times, we conclude that there exists a map $\sigma: \overline{1,m}\to \overline{1,n}$  such that for any
$i\in\overline{1,m}$, there is a unitary $u_i\in M$ satisfying
$u_i L(\Gamma_i) u_i^* \subset L(\Lambda_{\sigma(i)})$.

In a similar way, we obtain
a map $\tau: \overline{1,n}\to \overline{1,m}$ and a unitary $w_j\in M$, for any $j\in\overline{1,n}$, such that 
$w_jL(\Lambda_j) w_j^* \subset L(\Gamma_{\tau(j)})$, for any $j\in\overline{1,n}$. Thus, $u_{\tau(j)}w_j L(\Lambda_j) w_j^* u_{\tau(j)}^*\subset L(\Lambda_{\sigma(\tau(j))})$, for any $j\in\overline{1,n}$. By applying \cite[Theorem 1.2.1]{IPP05} we deduce that $\sigma\circ\tau ={\rm Id}$ and $u_{\tau(j)}w_j\in L(\Lambda_j)$, for any $j\in\overline{1,n}$. Similarly, we get $\tau\circ\sigma={\rm Id}$ and $w_{\sigma(i)}u_i\in L(\Gamma_i)$ for any $i\in\overline{1,m}$. In particular, $m=n$ and $u_i L(\Gamma_i) u_i^* = L(\Lambda_{\sigma(i)})$, for any $i\in\overline{1,n}$.
\hfill$\blacksquare$


    
    
    

\end{document}